\documentclass[reqno,11pt]{amsart}
\usepackage{multind}
\usepackage{hyperref}
\usepackage{graphicx}
\usepackage{amscd}
\usepackage{amssymb}
\usepackage[mathscr]{eucal}
\numberwithin{equation}{section}
\allowdisplaybreaks[1]

\title[Airy and Parabolic Cylinder Functions]
{Absence of Zeros and Asymptotic Error Estimates for Airy and Parabolic Cylinder Functions}

\author[F.\ Finster]{Felix Finster}
\thanks{F.F.\ is supported in part by the Deutsche Forschungsgemeinschaft.}
\address{Fakult\"at f\"ur Mathematik \\ Universit\"at Regensburg \\ D-93040 Regensburg \\ Germany}
\email{finster@ur.de}

\author[J.\ Smoller]{Joel Smoller \\ \\ October 2012}
\thanks{J.S.\ is supported in part by the National Science Foundation,
Grant No.\ DMS-110189.}
\address{Mathematics Department \\ The University of Michigan \\ Ann Arbor, MI 48109, USA}
\email{smoller@umich.edu}

\newtheorem{Def}{Def.}[section]
\newtheorem{Thm}[Def]{Theorem}
\newtheorem{Prp}[Def]{Proposition}
\newtheorem{Lemma}[Def]{Lemma}
\newtheorem{Remark}[Def]{Remark}
\newtheorem{Corollary}[Def]{Corollary}

\newcommand{\Thanks}{\vspace*{.5em} \noindent \thanks}
\newcommand{\Proof}{\begin{proof}}
\newcommand{\QED}{\end{proof} \noindent}
\newcommand{\QEDrem}{\ \hfill $\Diamond$}

\newcommand{\C}{\mathbb{C}}
\newcommand{\R}{\mathbb{R}}

\newcommand{\beq}{\begin{equation}}
\newcommand{\beql}[1]{\begin{equation} \label{#1}}
\newcommand{\eeq}{\end{equation}}

\newcommand{\WKB}{{\text{\tiny{\rm{WKB}}}}}

\newcommand{\A}{\mathfrak{A}}
\renewcommand{\O}{\mathscr{O}}

\DeclareMathOperator{\re}{Re}
\DeclareMathOperator{\im}{Im}
\DeclareMathOperator{\Ai}{Ai}
\DeclareMathOperator{\Bi}{Bi}

\begin{document}

\maketitle

\begin{abstract}
We derive WKB approximations for a class of Airy and parabolic cylinder functions in the complex plane,
including quantitative error bounds. We prove that all zeros of the Airy
function lie on a ray in the complex plane, and that the parabolic cylinder functions
have no zeros. We also analyze the Airy and Airy-WKB limit of the parabolic cylinder functions.
\end{abstract}

\tableofcontents

\section{Introduction and Statement of Results}
This paper represents an important first step in our program of proving the linearized stability of the Kerr
black hole.
The mathematical problem can be framed as showing that solutions of the so-called Teukolsky
equation~\cite{teukolsky, chandra} for spin s=2, with smooth compactly supported initial data outside the black hole, decay uniformly on compact sets.  The Teukolsky equation separates into an angular and radial ODE.  The angular equation involves a Sturm-Liouville operator of the form
\[ -\frac{d^2}{du^2} + W(u) \:. \]
Here the potential~$W$ is {\em{complex}}, and thus the operator is {\em{not}} self-adjoint.  Following the procedure developed in~\cite{angular}, we seek a spectral representation for this angular operator.
Our method for proving such a spectral representation requires detailed information on the eigensolutions.
This information can be obtained by ``glueing together'' approximate solutions
and controlling the error using the methods in~\cite{invariant}.
One method for obtaining such approximate solutions is to solve the Sturm-Liouville equation
\begin{equation} \label{5ode}
\Big( -\frac{d^2}{du^2} + V(u) \Big) \phi(u) = 0\:,
\end{equation}
where~$V$ is a linear or quadratic polynomial with complex coefficients.
The corresponding solutions are Airy and parabolic cylinder functions, respectively.
In this paper we analyze properties of these special functions.

Both of these special functions have interesting representations as contour integrals in the complex plane.
Since it is difficult to estimate these integral representations directly, we analyze them
with stationary phase-like methods
(for an introduction to these methods see for example~\cite[Section~7.7]{hormanderI})
to obtain approximate WKB solutions
\beq \label{WKBsol}
\frac{1}{\sqrt[4]{V}} \: \exp \Big( \int^u \sqrt{V} \Big)
\eeq
with quantitative error bounds.
In order to apply the methods in~\cite{invariant}, we also need to control the
function~$y:= \phi'/\phi$, which solves the associated Riccati equation
\beq  \label{riccati}
y' = V - y^2\:.
\eeq
For~$y$ to be well-behaved, we must show that~$\phi$ has no zeros.
This motivates our interest in ruling out zeros of the Airy and parabolic cylinder functions.

Apart from being of independent interest, the results obtained here will be used in the forthcoming
papers~\cite{tinvariant, tspectral}.

\section{Estimates for Airy Functions}
In this section, we assume that~$V$ is a linear function,
\[  V(u) = a + b u \qquad \text{with} \qquad a, b \in \C\:. \]
Then the Sturm-Liouville equation~\eqref{5ode} can be solved explicitly in terms of Airy functions,
\[ \phi(u) = \A(z) \qquad \text{where} \qquad
z = b^{-\frac{2}{3}} (a + b u) \]
and~$\A$ is a linear combination of the Airy functions~$\Ai$ and~$\Bi$ (see~\cite{DLMF}) 
\[ \A(z) := -i \sqrt{\pi} \,\big( \Ai(z) + i \,\Bi(z) \big) \]
(the reason for our specific linear combination is that it has a particularly simple WKB
asymptotic form; see~\eqref{AWKB} below).
The corresponding solution of the Riccati equation~\eqref{riccati} is given by
\[ y(u) = b^\frac{1}{3}\: \frac{\A'(z)}{\A(z)} \:. \]
Using the integral representations in~\cite[eqns~(9.5.4) and~(9.5.5)]{DLMF}
\begin{align}
\Ai(z) &= \frac{1}{2 \pi i} \int_{\infty e^{-\pi i/3}}^{\infty e^{\pi i/3}} \exp \left( \frac{t^3}{3} - z t \right) dt \label{Airep} \\
\Bi(z) &= \frac{1}{2 \pi} \int_{-\infty}^{\infty e^{\pi i/3}} \exp \left( \frac{t^3}{3} - z t \right) dt 
+ \frac{1}{2 \pi} \int_{-\infty}^{\infty e^{-\pi i/3}} \exp \left( \frac{t^3}{3} - z t \right) dt
\end{align}
(where~$\infty e^{i \varphi}$ refers to the end point of the contour~$\lim_{t \rightarrow +\infty} t e^{i \varphi}$,
etc.), we obtain
\beq \label{Arep}
\A(z) = \frac{1}{\sqrt{\pi}} \int_\Gamma \exp \left( \frac{t^3}{3} - z t \right) dt \:,
\eeq
where~$\Gamma$ is the contour~$\Gamma = -\R^+ \cup e^{-i \pi/3} \,\R^+$.
Note that the last integral is finite because the factor~$e^{t^3/3}$ decays exponentially
at both ends of the contour. From the computation
\begin{align*}
\A''(z) &= \frac{1}{\sqrt{\pi}} \int_\Gamma t^2 \:\exp \left( \frac{t^3}{3} - z t \right) dt \\
&= \frac{1}{\sqrt{\pi}} \int_\Gamma \left( z + \frac{d}{dt} \right) \exp \left( \frac{t^3}{3} - z t \right) dt
= z \,\A(z)
\end{align*}
one immediately verifies that~$\A$ is a solution of the Airy equation.

In the next lemma, we show that~$\A$ can be obtained from~$\Ai$ by a
rotation in the complex plane.
\begin{Lemma} \label{lemmaAAi} For any~$z \in \C$,
\[  \A \big( e^{-\frac{i \pi}{3}}\, z \big) = 2 \sqrt{\pi}\: e^{-\frac{i \pi}{6}}\:  \Ai(-z)\:. \]
\end{Lemma}
\Proof
We perform a change of variables in the integral~\eqref{Arep},
\begin{align*}
\A\big( e^{-\frac{i \pi}{3}}\, z \big) &=
\frac{1}{\sqrt{\pi}} \int_\Gamma \exp \left( \frac{t^3}{3} - e^{-\frac{i \pi}{3}}\,z t \right) dt
= \left\{ \tau = e^{\frac{2 \pi  i}{3}}\: t \right\} \\
&=  \frac{1}{\sqrt{\pi}} \:e^{-\frac{2 \pi i}{3}} \int_{\infty e^{-\pi i/3}}^{\infty e^{\pi i/3}}
\exp \left( e^{-2 \pi i} \:\frac{\tau^3}{3} - e^{-i \pi}\:z \tau \right) d\tau \\
&= \frac{1}{\sqrt{\pi}} \:e^{-\frac{2 \pi i}{3}} \int_{\infty e^{-\pi i/3}}^{\infty e^{\pi i/3}}
\exp \left( \frac{\tau^3}{3} + z \tau \right) d\tau \\
&\!\!\overset{\eqref{Airep}}{=} 2i \sqrt{\pi} \: e^{-\frac{2 \pi i}{3}}\: \Ai(-z)
= 2 \sqrt{\pi}\: e^{-\frac{i \pi}{6}}\: \Ai(-z) \:,
\end{align*}
giving the result.
\QED

\subsection{WKB Estimates} \label{secairyWKB}
\begin{figure}
\begin{picture}(0,0)%
\includegraphics{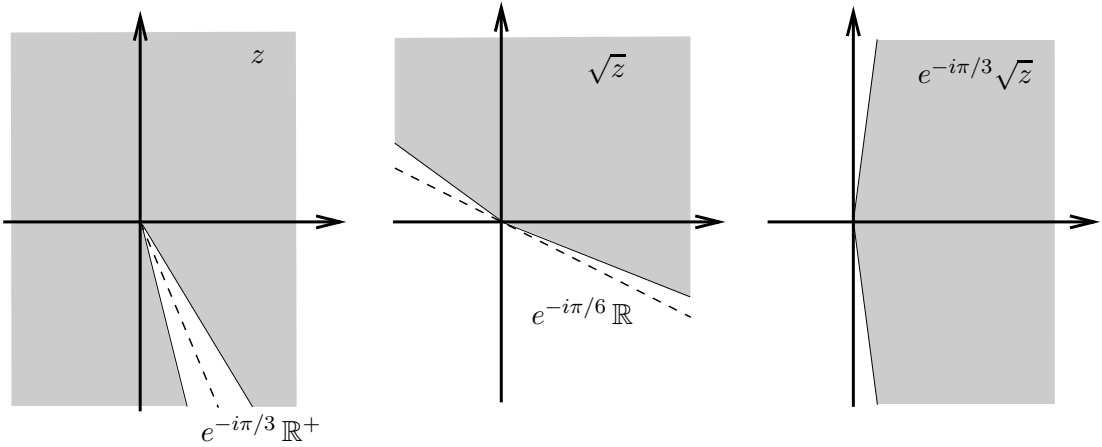}%
\end{picture}%
\setlength{\unitlength}{2486sp}%
\begingroup\makeatletter\ifx\SetFigFont\undefined%
\gdef\SetFigFont#1#2#3#4#5{%
  \reset@font\fontsize{#1}{#2pt}%
  \fontfamily{#3}\fontseries{#4}\fontshape{#5}%
  \selectfont}%
\fi\endgroup%
\begin{picture}(11006,4519)(-2792,-7742)
\put(-291,-3881){\makebox(0,0)[lb]{\smash{{\SetFigFont{11}{13.2}{\familydefault}{\mddefault}{\updefault}$z$}}}}
\put(-786,-7646){\makebox(0,0)[lb]{\smash{{\SetFigFont{11}{13.2}{\familydefault}{\mddefault}{\updefault}$e^{-i \pi/3} \,\R^+$}}}}
\put(2508,-6480){\makebox(0,0)[lb]{\smash{{\SetFigFont{11}{13.2}{\familydefault}{\mddefault}{\updefault}$e^{-i \pi/6} \,\R$}}}}
\put(3061,-4017){\makebox(0,0)[lb]{\smash{{\SetFigFont{11}{13.2}{\familydefault}{\mddefault}{\updefault}$\sqrt{z}$}}}}
\put(6363,-4078){\makebox(0,0)[lb]{\smash{{\SetFigFont{11}{13.2}{\familydefault}{\mddefault}{\updefault}$e^{-i \pi/3} \sqrt{z}$}}}}
\end{picture}%
\caption{Parameter regions in the complex plane.}
\label{figzbranch}
\end{figure}

Our first goal is to get asymptotic expansions and rigorous estimates of
the functions~$\A'(z)$ and~$\A(z)$.
We expect that for large~$z$, the Airy solution should go over to the
WKB wave function~\eqref{WKBsol} corresponding to the Airy potential~$V(z) = z$,
\beq \label{AWKB}
\A_\WKB(z) := z^{-\frac{1}{4}}\: e^{\frac{2}{3}\, z^\frac{3}{2}}\:.
\eeq
Here we define the roots by~$z^\alpha = \exp(\alpha \log(z))$, where
the logarithm has a branch cut along the ray
\beq \label{cut}
e^{-\frac{i \pi}{3}}\: \R^+ \:.
\eeq
In the next theorem, we show that with this branch cut, the WKB wave function~\eqref{AWKB}
approximates the Airy function~$\A$ for large~$z$, with rigorous error bounds.
\begin{Thm} \label{thmairy}
Assume that for given~$\varepsilon \in (0, \frac{\pi}{6})$,
\beq \label{zreim}
\arg z \not \in \Big(-\frac{i \pi}{3}-\varepsilon, -\frac{i \pi}{3}+\varepsilon \Big)  \mod 2 \pi \:.
\eeq
Then
\begin{align}
\left| \frac{\A(z)}{\A_\WKB(z)} - 1 \right| &\leq
\frac{3}{|z|^\frac{3}{4}}\: \frac{1}{\sin^2(\varepsilon/2)} \left(1 + \log^2 \big(1 +
\sin^\frac{3}{2}(\varepsilon/2)\, |z|^\frac{3}{4} \big) \right) \label{Aes} \\
\left| \frac{d}{dz} \left( \frac{\A(z)}{\A_\WKB(z)} \right) \right| &\leq
\frac{2}{|z|^\frac{7}{4}}\: \frac{1}{\sin^3(\varepsilon/2)} \left(1 + \log^3 \big(1 +
\sin^\frac{3}{2}(\varepsilon/2)\, |z|^\frac{3}{4} \big) \right) . \label{Apes}
\end{align}
\end{Thm}
\Proof The assumption~\eqref{zreim} and our branch convention for the square root imply that the
parameter~$c$ defined by
\beq \label{cdef}
c = \re \left( e^{-\frac{i \pi}{3}}\, \sqrt{z} \right)
\eeq
is positive and bounded by (see Figure~\ref{figzbranch})
\beq \label{cbound}
c \geq |z|^\frac{1}{2}\, \sin (\varepsilon/2 \big)\:.
\eeq

Rewriting the argument of the exponential in~\eqref{Arep} as
\[ \frac{t^3}{3} - z t = \frac{2}{3}\: z^\frac{3}{2} - \sqrt{z}\: \tau^2 + \frac{\tau^3}{3} 
\qquad \text{with} \qquad \tau := t + \sqrt{z}\:, \]
we obtain
\begin{align}
\A(z) &= \frac{1}{\sqrt{\pi}} \int_{e^{-i \pi/6}\,\R}
\exp \left( \frac{2}{3}\: z^\frac{3}{2} - \sqrt{z}\: \tau^2 + \frac{\tau^3}{3} \right) d\tau \label{A1} \\
&= \frac{1}{\sqrt{\pi}} \: e^{-\frac{i \pi}{6}} \int_{-\infty}^\infty
\exp \left( \frac{2}{3}\: z^\frac{3}{2} - e^{-\frac{i \pi}{3}} \sqrt{z}\: s^2 - \frac{i s^3}{3} \right) ds \:, \label{A2}
\end{align}
where in~\eqref{A1} we deformed the contour, whereas in the last line we
introduced the new integration variable~$s = e^{i \pi/6}\, \tau$.
Note that in~\eqref{A2} the factor~$e^{-i s^3/3}$ is merely a phase,
but the quadratic term in the exponent still ensures convergence of the integral in view of~\eqref{cdef}.

If we drop the cubic term in the exponent, the resulting Gaussian integral
can be computed to give precisely the WKB wave function,
\[ \A_\WKB(z) = \frac{1}{\sqrt{\pi}} \: e^{-\frac{i \pi}{6}} \int_{-\infty}^\infty \exp \left( \frac{2}{3}\: z^\frac{3}{2} - e^{-\frac{i \pi}{3}}\, \sqrt{z}\, s^2 \right) ds \:. \]
We thus obtain the error term
\beq \begin{split}
E(z) &:= \A(z) - \A_\WKB(z) \\
&\,= \frac{1}{\sqrt{\pi}} \: e^{-\frac{i \pi}{6}} \int_{-\infty}^\infty
\exp \left( \frac{2}{3}\: z^\frac{3}{2} - e^{-i \pi/3} \sqrt{z}\: s^2 \right)
\left( e^{-\frac{i s^3}{3}} - 1 \right) ds\:.
\end{split} \label{Ezdef}
\eeq
After estimating the error term by
\[ \left| \sqrt{\pi} \:e^{-\frac{2}{3}\, z^\frac{3}{2}}\: E(z) \right| \leq
\int_{-\infty}^\infty e^{-c s^2} \left| e^{-\frac{i s^3}{3}} - 1 \right| ds \:, \]
we decompose the integral into integrals over the regions~$[-L,L]$
and~$\R \setminus [-L,L]$. We bound these integrals as follows,
\begin{align}
\int_{\R \setminus [-L,L]} e^{-c s^2} \left| e^{-\frac{i s^3}{3}} - 1 \right| ds
&\leq 2 \int_{\R \setminus [-L,L]} e^{-c s^2}\: ds \nonumber \\
&\leq 2 \,e^{-\frac{c L^2}{2}} \int_{-\infty}^\infty e^{-\frac{c s^2}{2}}\: ds
= \frac{2 \sqrt{2 \pi}}{\sqrt{c}}\,e^{-\frac{c L^2}{2}} \label{ies1} \\
\int_{-L}^L e^{-c s^2} \left| e^{-\frac{i s^3}{3}} - 1 \right| ds
&\leq 2 \int_0^L \frac{s^3}{3}\: ds \leq
\frac{L^4}{6}\:, \label{ies2}
\end{align}
where in the last line we used that for all~$x \in \R$ the inequality~$|e^{ix} -1| \leq |x|$
holds. In order to make both errors of about the same size, we choose
\beq \label{Lchoice}
L = \sqrt{\frac{2}{c}}\: \log^\frac{1}{2} \big( 1 + c^\frac{3}{2} \big)\:.
\eeq
We thus obtain
\[ \left| e^{-\frac{2}{3}\, z^\frac{3}{2}}\: E(z) \right| \leq
\frac{2^\frac{3}{2}}{c^2}
+ \frac{2}{3 \,\sqrt{\pi}}\: \frac{1}{c^2}\:\log^2 \big(1 + c^{3/2} \big)
\leq \frac{3}{c^2} \left( 1 + \log^2 \big(1 + c^{3/2} \big) \right) . \]
Applying~\eqref{cbound} and~\eqref{AWKB} gives~\eqref{Aes}.

In order to derive~\eqref{Apes}, we introduce the abbreviations
\begin{align}
I(z) &=  \int_{-\infty}^\infty \exp \left(- e^{-\frac{i \pi}{3}} \sqrt{z}\: s^2 - \frac{i s^3}{3} \right) ds \label{IA} \\
I_\WKB(z) &= \int_{-\infty}^\infty \exp \left(- e^{-\frac{i \pi}{3}}\, \sqrt{z}\, s^2 \right) ds 
= e^{\frac{i \pi}{6}}\: \sqrt{\pi}\: z^{-\frac{1}{4}} \:. \label{Igauss}
\end{align}
Then
\[ \frac{d}{dz} \left( \frac{\A(z)}{\A_\WKB(z)} \right)
= \frac{d}{dz} \left( \frac{I(z)}{I_\WKB(z)} \right)
= \frac{I'(z) - I'_\WKB(z)}{I_\WKB(z)} + \frac{I_\WKB'(z)}{I_\WKB(z)}\:
\frac{\A_\WKB(z) - \A(z)}{\A_\WKB(z)} \:. \]
Using~\eqref{Igauss} and~\eqref{Aes}, the second summand is estimated by
\beql{ques1}
\left| \frac{I_\WKB'(z)}{I_\WKB(z)}\:
\frac{\A_\WKB(z) - \A(z)}{\A_\WKB(z)} \right| \leq
\frac{1}{|z|^\frac{7}{4}}\: \frac{1}{\sin^2(\varepsilon/2)} \left(1 + \log^2 \big(1 +
\sin^\frac{3}{2}(\varepsilon/2) |z|^\frac{3}{4} \big) \right) .
\eeq
In order to estimate the first summand, we differentiate~\eqref{IA} and~\eqref{Igauss}
to obtain, similar to~\eqref{Ezdef}, the identity
\begin{align*}
\sqrt{z}& \left( I'(z) - I_\WKB'(z) \right)
= \int_{-\infty}^\infty
\left(- \frac{1}{2} \:e^{-\frac{i \pi}{3}} \, s^2 \right) \exp \left(- e^{-\frac{i \pi}{3}} \sqrt{z}\: s^2 \right)
\left( e^{-\frac{i s^3}{3}} - 1 \right) ds\:.
\end{align*}
It follows that
\[ \left| \sqrt{z} \,\big( I'(z) - I_\WKB'(z) \big) \right|
\leq \frac{1}{2} \int_{-\infty}^\infty s^2 \:e^{-c s^2} \left| e^{-\frac{i s^3}{3}} - 1 \right| ds\:. \]
Estimating the integral exactly as in~\eqref{ies1} and~\eqref{ies2}, we obtain
\[ \frac{1}{2} \int_{-\infty}^\infty s^2 \:e^{-c s^2} \left| e^{-\frac{i s^3}{3}} - 1 \right| ds \leq
 \frac{\sqrt{\pi}}{c^\frac{3}{2}}\,e^{-\frac{c L^2}{2}} + \frac{L^6}{18} \:. \]
Choosing~$L$ again according to~\eqref{Lchoice}, we obtain
\begin{align*}
\left| \sqrt{z} \,\big( I'(z) - I_\WKB'(z) \big) \right|
&\leq \frac{1}{\sqrt{\pi}\: c^3}  \left( 1
+ \log^3 (1+ c^\frac{3}{2}) \right) \:.
\end{align*}
Applying~\eqref{Igauss} and~\eqref{cbound}, we conclude that
\[ \left| \frac{I'(z) - I'_\WKB(z)}{I_\WKB(z)} \right|
\leq \frac{1}{|z|^\frac{7}{4}}\: \frac{1}{\sin^3(\varepsilon/2)} \left(1 + \log^3 \big(1 +
\sin^\frac{3}{2}(\varepsilon/2)\, |z|^\frac{3}{4} \big) \right) . \]
Noting that this contribution dominates~\eqref{ques1}, we obtain~\eqref{Apes}.
\QED

In particular, this theorem allows us to take the limit as~$z$ goes to infinity
along a ray through the origin, as long as we stay away from the branch cut~\eqref{cut}.
\begin{Corollary} \label{cor22}
Let~$z_0 \in \C$ be off the ray~$e^{-\frac{i \pi}{3}} \R^+$. Then
\[ \lim_{t \rightarrow \infty} t^{-\frac{1}{2}} \,\frac{\A'(t z_0)}{\A(t z_0)} = \sqrt{z_0}\:, \]
where we again used the sign convention for the square root introduced just before~\eqref{cut}.
\end{Corollary}
\Proof We choose~$\varepsilon$ so small that~$z_0$ satisfies the condition~\eqref{zreim}. Then
\[ \lim_{t \rightarrow \infty} t^{-\frac{1}{2}} \,\frac{\A'(t z_0)}{\A(t z_0)}
= \lim_{t \rightarrow \infty} t^{-\frac{1}{2}} \,\frac{\A_\WKB'(t z_0)}{\A_\WKB(t z_0)} \:, \]
and computing~$\A'_\WKB$ from~\eqref{AWKB} gives the result.
\QED

\subsection{An Estimate for the Riccati Equation for a Complex Potential} \label{sec22}
Our next goal is to show that the Airy function~$\A$ has no zeros except on the branch cut~\eqref{cut}.
In preparation for this, we now derive an estimate for solutions of the Riccati equation
for a general potential~$V$ with~$\im V \geq 0$.
Thus let~$y$ be a solution of the Riccati equation~\eqref{riccati}.
Decomposing~$y$ into its real and imaginary parts, $y=\alpha+i \beta$, we obtain the system
\begin{align}
\alpha' &= \re V - \alpha^2 + \beta^2 \label{aeq} \\
\beta' &= \im V - 2 \alpha \beta \:. \label{beq}
\end{align}
Moreover, we set
\beq \label{sigmadef}
\sigma = \exp \left( 2  \int^u \alpha \right)\:.
\eeq

Let us assume that~$\im V \geq 0$ on an interval~$[u_0,u_1]$ and that~$(\alpha, \beta)$ is
a solution on a closed subinterval~$[u_0, u_2] \subset [u_0, u_1]$.
If~$\beta(u_0)>0$, it follows immediately from~\eqref{beq} that~$\beta > 0$ on the whole interval~$[u_0, u_2]$.
A short calculation yields
\[ (\sigma \beta)'(u) = \sigma \im V \geq 0 \:. \]
Hence the function~$\sigma \beta$ is monotone increasing and
\[ \beta(u) \geq \beta_0 \,\sigma_0\: \frac{1}{\sigma(u)}\:, \]
where~$\beta_0 = \beta(u_0)$ and similarly for all other functions.
Using this inequality in~\eqref{aeq}, we obtain
\beq \label{aineq}
\alpha' \geq \re V - \alpha^2 + \frac{\beta_0^2 \,\sigma_0^2}{\sigma(u)^2}\:.
\eeq
Integrating this inequality gives the following estimate.
\begin{Lemma} Assume that~$\im V \geq 0$ on the interval~$[u_0, u_1]$
and that~$(\alpha, \beta)$ is a solution on a closed subinterval~$[u_0, u_2] \subset [u_0, u_1]$
with~$\beta(u_0)>0$.
Then the function~$\alpha$ satisfies the upper bound
\[ \alpha(u) \leq g(u) \qquad \text{on~$[u_0, u_2]$}\:, \]
where we set
\[ g(u) = \frac{\sqrt{C} \sinh(x) + \alpha(u_2)\, \cosh(x)}{\cosh(x) + \alpha(u_2)\, \sinh(x) / \sqrt{C}} \]
and
\[ C = \inf_{[u_0, u_1]} \re V \in \R \:,\qquad x = \sqrt{C}\, (u-u_1) \:. \]
\end{Lemma}
\Proof Dropping the last summand in~\eqref{aineq}, we obtain the inequality
\[ \alpha' \geq C - \alpha^2  \:. \]
Then~$\alpha \leq g$, where~$g$ is the solution of the corresponding differential equation
\[ g' = C - g^2\:,\qquad g(u_2) = \alpha(u_2)\:. \]
Solving this differential equation gives the result.
\QED
Using the result of this lemma in~\eqref{sigmadef}, we get
\begin{align*}
\sigma(u) = \sigma_0\, \exp \left( 2 \int_{u_0}^u \alpha \right) \leq
\sigma_0\, \exp \left( 2 \int_{u_0}^u g \right) .
\end{align*}
Computing the integral, we obtain
\[ \sigma(u) \leq \sigma_0 \: \left( \frac{\sqrt{C} \cosh \big( \sqrt{C} (u-u_2) \big) + \alpha(u_2)\, \sinh \big(\sqrt{C} (u-u_2) \big)}
{\sqrt{C} \cosh \big(\sqrt{C} (u-u_0) \big) + \alpha(u_2)\, \sinh \big( \sqrt{C} (u-u_0) \big)}
\right)^2 . \]
Using this estimate in~\eqref{aineq} and setting~$u=u_2$ gives the following result.

\begin{Prp} \label{prp23} Assume that~$\im V \geq 0$ on~$[u_0, u_1]$ 
and that~$\beta(u_0)>0$. Assume furthermore that the solution~$(\alpha, \beta)$
exists on the interval~$[u_0, u] \subset [u_0, u_1]$ and that~$\alpha'(u) \leq 0$. Then
\[ 0 \geq C - \alpha(u)^2 + \beta_0^2
\left[ \cosh \big( \sqrt{C} (u-u_0) \big) + \frac{\alpha(u)}{\sqrt{C}}\: \sinh \big(\sqrt{C} (u-u_0) \big)
\right]^4 \:. \]
\end{Prp}

\begin{Corollary} \label{cor25}
Assume that either~$\im V|_{[u_0, u_1]} \geq 0$ and~$\beta(u_0)>0$
or~$\im V|_{[u_0, u_1]} \leq 0$ and~$\beta(u_0)<0$.
Then the solution~$(\alpha, \beta)$ exists and is bounded on the interval~$[u_0, u_1]$.
\end{Corollary}
\Proof It suffices to consider the case~$\im V|_{[u_0, u_1]} \geq 0$ and~$\beta(u_0)>0$
because the other case is obtained by taking the complex conjugate of the Riccati equation~\eqref{riccati}.

Let us assume conversely that the solution~$(\alpha, \beta)$ blows up at a point~$u_2 \in [u_0, u_1]$.
We can assume that~$\alpha$ blows up, because otherwise~$\beta$ could be obtained by
integrating~\eqref{beq},
\beq \label{varc}\beta(u) =  e^{-2 \int_{u_0}^u \alpha}\, \beta(u_0)  + \int_{u_0}^u e^{-2 \int_\tau^u \alpha}
\im V(\tau)\: d\tau \:.
\eeq
Assume that~$\alpha$ tends to $+\infty$. Then~\eqref{varc}  shows that~$\beta$ stays bounded,
and thus the right hand side of~\eqref{aeq} tends to~$-\infty$, a contradiction.
On the other hand, if~$\alpha$ is unbounded from below, then there is a sequence~$v_n \in [u_0, u_2)$
with~$\alpha(v_n) \rightarrow -\infty$ and~$\alpha'(v_n) \leq 0$.
This contradicts Proposition~\ref{prp23}.
\QED

\subsection{Locating the Zeros of the Airy Function} \label{secairyzero}
Our method for ruling out zeros of the Airy function~$\A$  at a given point~$z_0 \in \C$
is to consider the Airy function along straight lines of the form
\beq \label{zudef}
z(u) = z_0 + \lambda u \qquad \text{with} \qquad u \in \R\:,
\eeq
where~$\lambda$ is a complex parameter. Then the function~$\phi(u) := \A(z(u))$
satisfies the Sturm-Liouville equation~\eqref{5ode} with
\[ V(u) = \lambda^2\; (z_0 + \lambda u) \:. \]
Applying Corollary~\ref{cor25}
to the corresponding Riccati solution~$y:=\phi' / \phi$
yields the following proposition. This proposition follows immediately from Lemma~\ref{lemmaAAi}
and the fact that the function~$\Ai$ only has zeros on the negative real axis
(see~\cite[\S9.9(i)]{DLMF}). Nevertheless, we present a proof
in order to illustrate the methods of Section~\ref{sec22} in a simple example
(these methods will be used again in Section~\ref{sec32}).

\begin{Prp} \label{thmairy2} The function~$\A$ has no zeros off the ray~$e^{-i \pi/3} \,\R^+$.
\end{Prp}
\Proof Let~$z_0$ be a point off the ray~$e^{-i \pi/3} \,\R^+$.
We choose~$\lambda$ such that~$(-\lambda)$ is also off the ray~$e^{-i \pi/3} \,\R^+$.
Then Corollary~\ref{cor22} applies and yields
\[ \lim_{u \rightarrow -\infty} |u|^{-\frac{1}{2}} \,y(u) 
= \lambda \lim_{u \rightarrow -\infty} |u|^{-\frac{1}{2}} \,\frac{\A'(z(u))}{\A(z(u))}
= \lambda\:\sqrt{-\lambda} \:. \]
Using our sign convention for the square root, we may choose~$\lambda$ such that
\beq \label{phirange}
\sqrt{-\lambda} = e^{i \varphi} \qquad \text{with} \qquad \varphi \in \Big(\frac{\pi}{6}, \frac{\pi}{3} \Big)
\cup \Big( \frac{\pi}{3}, \frac{\pi}{2} \Big) \:;
\eeq
then
\beq \label{limy}
\lim_{u \rightarrow -\infty} |u|^{-\frac{1}{2}} \,y(u) = -e^{3 i \varphi}
\eeq
(note that the parameter~$-\lambda = e^{2 i \varphi}$ is indeed off the ray~$e^{-i \pi/3} \,\R^+$
because~$\pi/3$ is excluded in~\eqref{phirange}).
Moreover,
\beq \label{imV}
\im V(u) = \im ( \lambda^2 z_0)  + u \im ( \lambda^3) 
= \im ( e^{4 i \varphi}  z_0)  - u \im ( e^{6 i \varphi}) \:.
\eeq

\begin{figure}
\begin{picture}(0,0)%
\includegraphics{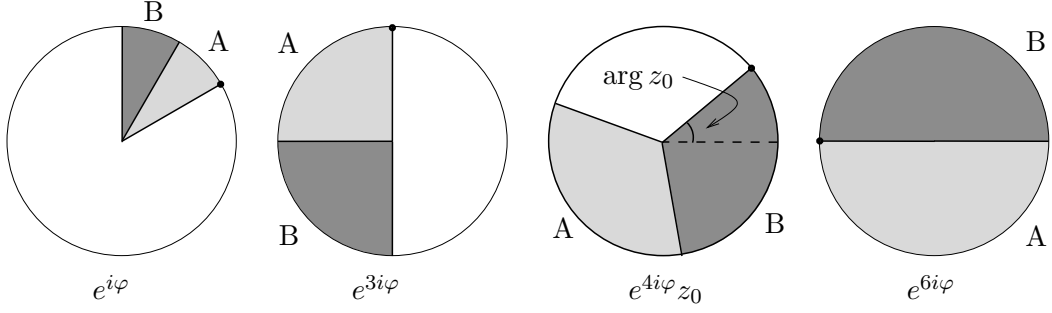}%
\end{picture}%
\setlength{\unitlength}{2486sp}%
\begingroup\makeatletter\ifx\SetFigFont\undefined%
\gdef\SetFigFont#1#2#3#4#5{%
  \reset@font\fontsize{#1}{#2pt}%
  \fontfamily{#3}\fontseries{#4}\fontshape{#5}%
  \selectfont}%
\fi\endgroup%
\begin{picture}(10416,3085)(-3400,-7132)
\put(-1391,-4561){\makebox(0,0)[lb]{\smash{{\SetFigFont{11}{13.2}{\familydefault}{\mddefault}{\updefault}A}}}}
\put(-2032,-4257){\makebox(0,0)[lb]{\smash{{\SetFigFont{11}{13.2}{\familydefault}{\mddefault}{\updefault}B}}}}
\put(6765,-6510){\makebox(0,0)[lb]{\smash{{\SetFigFont{11}{13.2}{\familydefault}{\mddefault}{\updefault}A}}}}
\put(6768,-4567){\makebox(0,0)[lb]{\smash{{\SetFigFont{11}{13.2}{\familydefault}{\mddefault}{\updefault}B}}}}
\put(-683,-6494){\makebox(0,0)[lb]{\smash{{\SetFigFont{11}{13.2}{\familydefault}{\mddefault}{\updefault}B}}}}
\put(-695,-4583){\makebox(0,0)[lb]{\smash{{\SetFigFont{11}{13.2}{\familydefault}{\mddefault}{\updefault}A}}}}
\put(5581,-7036){\makebox(0,0)[lb]{\smash{{\SetFigFont{11}{13.2}{\familydefault}{\mddefault}{\updefault}$e^{6 i \varphi}$}}}}
\put( 45,-7036){\makebox(0,0)[lb]{\smash{{\SetFigFont{11}{13.2}{\familydefault}{\mddefault}{\updefault}$e^{3 i \varphi}$}}}}
\put(-2538,-7036){\makebox(0,0)[lb]{\smash{{\SetFigFont{11}{13.2}{\familydefault}{\mddefault}{\updefault}$e^{i \varphi}$}}}}
\put(2806,-7036){\makebox(0,0)[lb]{\smash{{\SetFigFont{11}{13.2}{\familydefault}{\mddefault}{\updefault}$e^{4 i \varphi} z_0$}}}}
\put(4164,-6381){\makebox(0,0)[lb]{\smash{{\SetFigFont{11}{13.2}{\familydefault}{\mddefault}{\updefault}B}}}}
\put(2046,-6431){\makebox(0,0)[lb]{\smash{{\SetFigFont{11}{13.2}{\familydefault}{\mddefault}{\updefault}A}}}}
\put(2548,-4886){\makebox(0,0)[lb]{\smash{{\SetFigFont{11}{13.2}{\familydefault}{\mddefault}{\updefault}$\arg z_0$}}}}
\end{picture}%
\caption{Admissible phases.}
\label{figcircle}
\end{figure}

In the considered range of~$\varphi$, the imaginary parts of~$e^{3 i \varphi}$ and~$e^{6 i \varphi}$
have opposite signs (see Figure~\ref{figcircle}).
As a consequence, 
one sees from~\eqref{limy} and~\eqref{imV} that for sufficiently large negative~$u$,
the functions~$\im y$ and~$\im V$ are non-zero and have the same sign.
Moreover, we can choose~$\varphi$ such that the imaginary parts of~$e^{4 i \varphi}  z_0$
and~$e^{6 i \varphi}$ have the same sign (see Figure~\ref{figcircle}).
We conclude from~\eqref{imV} that the imaginary part of~$V$
has fixed sign on the interval~$(\infty, 0]$.
Thus Corollary~\ref{cor25} shows that~$y(u)$ is bounded on the interval~$(-\infty, 0]$.
Hence the corresponding Sturm-Liouville solution~$\phi(u) = \exp(\int^u y)$ has no zeros
on this interval. In particular, we conclude that~$0 \neq \phi(0) = \A(z_0)$.
\QED

\section{Estimates for Parabolic Cylinder Functions}
We now consider a quadratic potential
\[ V(u) = \alpha + \frac{\beta}{4}\: (u - \gamma)^2 \qquad \text{with} \qquad \alpha, \beta, \gamma \in \C\:. \]
The corresponding differential equation~\eqref{5ode} can be solved explicitly in terms
of the parabolic cylinder function, as we now recall. The parabolic cylinder function,
which we denote by~$U_a(z)$, is a solution of the differential equation
\beq \label{parabol}
U_a''(z) = \Big( \frac{z^2}{4}\: +a \Big) \, U_a(z)\:.
\eeq
Setting
\[ \phi(u) = U_a(z) \qquad \text{with} \qquad
a = \frac{\alpha}{\sqrt{\beta}} \:,\quad z = \beta^{\frac{1}{4}}\, (u-\gamma) \:, \]
a short calculation shows that~$\phi$ indeed satisfies~\eqref{5ode}.

According to~\cite[\S12.5(ii)]{DLMF}, particular solutions of~\eqref{parabol} can be written as contour integrals,
\beq \label{Uapm}
U_a^\pm(z) = e^{-\frac{z^2}{4}} \int_{\Gamma_\pm} e^{z t-\frac{t^2}{2}}\: t^{a-\frac{1}{2}}\: dt\:.
\eeq
We consider the special solution obtained by choosing
\[ \Gamma_\pm = \R \pm i \]
as the contour from~$-\infty \pm i$ to~$\infty \pm i$,
where we take the usual convention~$t^{a-\frac{1}{2}}
= \exp((a-\frac{1}{2}) \log t)$, and the logarithm has its branch cut along the negative real axis.
Whenever we omit the indices~$\pm$, our arguments apply to both cases.
It is verified by a direct computation that~$U_a$ satisfies~\eqref{parabol}; namely,
\begin{align*}
U_a'(z) &= \int_\Gamma \left(-\frac{z}{2}+t \right) e^{-\frac{z^2}{4}+z t-\frac{t^2}{2}}\: t^{a-\frac{1}{2}}\: dt \\
U_a''(z) &= \int_\Gamma \left(\frac{z^2}{4}-\frac{1}{2} -t\, (z-t) \right) e^{-\frac{z^2}{4}+z t-\frac{t^2}{2}}\: t^{a-\frac{1}{2}}\: dt \\
&= \int_\Gamma \: t^{a-\frac{1}{2}} \left(\frac{z^2}{4}-\frac{1}{2} -t\, \frac{\partial}{\partial t} \right)
e^{-\frac{z^2}{4}+z t-\frac{t^2}{2}}\: dt \\
&= \int_\Gamma \left(\frac{z^2}{4}-\frac{1}{2} + \Big(a+\frac{1}{2} \Big)  \right)\: t^{a-\frac{1}{2}} 
e^{-\frac{z^2}{4}+z t-\frac{t^2}{2}} \: dt = \Big( \frac{z^2}{4}\: +a \Big) \,U_a(z) \:,
\end{align*}
where in the last line we integrated by parts.

In the special cases
\beql{aexcept}
a = \frac{1}{2}\:, \frac{3}{2}, \frac{5}{2}, \ldots\:,
\eeq
the integral~\eqref{Uapm} is Gaussian and can easily be computed in closed form
(giving the well-known eigensolutions of the harmonic oscillator).
Therefore, we may restrict attention to the case when~\eqref{aexcept} is violated.
As~\eqref{parabol} only involves~$z^2$, it is obvious that~$U_a(-z)$ is also a parabolic cylinder
function. Indeed, by a change of variables we see that
\begin{align*}
U_a^+(-z) &= e^{-\frac{z^2}{4}} \int_{\Gamma_+} e^{-z t-\frac{t^2}{2}}\: t^{a-\frac{1}{2}}\: dt
= \left\{ \tau = - t \right\} \\
&= e^{-\frac{z^2}{4}} \int_{\Gamma_-} e^{z \tau-\frac{\tau^2}{2}}\: (-\tau)^{a-\frac{1}{2}}\: d\tau \:.
\end{align*}
Using that~$\log(-\tau) = \log(\tau) - i \pi$, we obtain the relation
\beql{Uapmrel}
U_a^+(-z) = e^{-i \pi (a-\frac{1}{2})}\: U_a^-(z) \:.
\eeq
Moreover, an elementary computation shows that
\beql{Uaprel} \begin{split}
U_a^+(0) &= 2^{a-\frac{1}{2}}\: (1-i e^{i \pi a})\: \Gamma \Big( \frac{a}{2} + \frac{1}{4} \Big) \\
(U_a^+)'(0) &= 2^{a+\frac{1}{2}}\: (1+i e^{i \pi a})\: \Gamma \Big( \frac{a}{2} + \frac{3}{4} \Big) \:.
\end{split}
\eeq
Thus the Wronskian of~$U_a^+$ and~$U_a^-$ is given by
\begin{align*}
w(U_a^+, U_a^-) &\overset{\eqref{Uapmrel}}{=} 2 e^{i \pi (a-\frac{1}{2})} (U_a^+)'(0)\: U_a^+(0) \\
&\;= e^{i \pi (a-\frac{1}{2})}\: 2^{2 a +1} \:(1+e^{2 i \pi a})\: \Gamma \Big( \frac{a}{2} + \frac{1}{4} \Big)
\:\Gamma \Big( \frac{a}{2} + \frac{3}{4} \Big) \:.
\end{align*}
From these formulas, one sees that~$U_a^+$ and~$U_a^-$ are linearly independent
except in the trivial cases~\eqref{aexcept}
(note that for the values~$a=-\frac{1}{2}, -\frac{3}{2}, \ldots$, the poles of the~$\Gamma$ functions
in~\eqref{Uaprel} are cancelled by the zeros of the corresponding factors~$(1 \pm i e^{i \pi a})$).

\subsection{WKB Estimates} \label{secparWKB}
Integral representations for parabolic cylinder function have been used
to derive asymptotic expansions (see for example~\cite{olver} and~\cite{temme+vidunas}).
The goal of this section is to derive rigorous estimates with quantitative error bounds
which go beyond~\cite{olver, temme+vidunas} and
show that in a certain parameter range, the function~$U_a(z)$ is well-approximated by a WKB wave function.
More precisely, the WKB wave function~\eqref{WKBsol} becomes
\beq \begin{split}
\phi_\WKB(z) &= \left( \frac{4}{z^2+4a} \right)^{\frac{1}{4}} \exp \left( \frac{1}{2} \int^z \sqrt{ \zeta^2 + 4a  }\;
d\zeta \right) \\
&= \frac{\big( z + \sqrt{z^2-b+2}\, \big)^{\frac{1}{2}-\frac{b}{4}}}{(z^2-b+2)^\frac{1}{4}}\:
\exp \left( {\frac{z}{4}\: \sqrt{z^2-b+2}} \right) ,
\end{split} \label{parWKB}
\eeq
where we introduced the abbreviation
\[ b=-4 \,\Big(a-\frac{1}{2} \Big) \:. \]

We write the parabolic cylinder function as
\beq
U_a(z) = e^{-\frac{z^2}{4}} \int_\Gamma e^{f(t)}\: dt \:, \label{Uaz}
\eeq
where the function~$f$ and its derivatives are given by
\begin{align}
f(t) &= z t - \frac{t^2}{2} -\frac{b}{4}\: \log t\:, &
f'(t) &= z - t - \frac{b}{4t} \label{f01} \\
f''(t) &= -1 + \frac{b}{4 t^2} \:, &
f'''(t) &= -\frac{b}{2 t^3} \:. \label{f23}
\end{align}
The zeros of~$f'$ are computed to be
\beq \label{tpm}
t_\pm = \frac{1}{2} \left(z \pm \sqrt{z^2 -b} \right) .
\eeq
We choose~$t_0$ equal to either~$t_+$ or~$t_-$.
As we excluded the special cases~\eqref{aexcept}, we know that~$b \neq 0$
and thus~$t_\pm \neq 0$.

In order to obtain a ``stationary phase-type'' approximation to~$U_a$,
we let~$\tilde{f}$ be the quadratic Taylor approximation of~$f$,
\begin{eqnarray*}
\tilde{f}(t) = f(t_0) + \frac{f''(t_0)}{2}\: (t-t_0)^2 \:.
\end{eqnarray*}
It is convenient to introduce the parametrization
\beq \label{tpar}
t = t_0 \,(1+\tau) \:;
\eeq
then
\begin{align}
f(t) - f(t_0) &= -\frac{t_0^2}{2}\: \tau^2 + \frac{b}{4} \,\big( \tau - \log(1+\tau) \big) \label{ft0rel} \\
\tilde{f}(t) - f(t_0) &= \left(-\frac{t_0^2}{2} + \frac{b}{8} \right) \tau^2 \:. \label{ftt0rel}
\end{align}
The coefficient of~$\tau^2$ in~\eqref{ft0rel} is given by
\beq \label{cdef2}
d := -\frac{t_0^2}{2} = -\frac{1}{8} \left(z \pm \sqrt{z^2-b} \right)^2 \:.
\eeq

We now deform the contour~$\Gamma$ to the the straight line
\beq \label{alphadef}
\tau(u) = e^{i \alpha} \,u \qquad \text{with~$u \in \R$ and~$\alpha \in [0, \pi) $}\:.
\eeq
Then
\begin{align*}
g(u) \;&\!\!:= \tau(u) - \log \big(1+\tau(u) \big) = e^{i \alpha} u - \log(1+e^{i \alpha} u ) \\
&= \int_0^u \frac{d}{ds} \Big( e^{i \alpha} s - \log(1+e^{i \alpha} s) \Big)\, ds
= \int_0^u \Big( e^{i \alpha} - \frac{e^{i \alpha}}{1+e^{i \alpha} s} \Big)\, ds \\
&= \int_0^u \frac{e^{2 i \alpha} s}{1+e^{i \alpha} s} \: ds
= \int_0^u \frac{e^{2 i \alpha} s\:(1+e^{-i \alpha} s)}{|1+e^{i \alpha} s|^2} \: ds
= \int_0^u \frac{e^{2 i \alpha} s+e^{i \alpha} s^2}{1+s^2 + 2 s \cos \alpha} \: ds \\
&= e^{2 i \alpha} I_1(u) + e^{i \alpha} I_2(u) \:,
\end{align*}
where~$I_1$ and~$I_2$ are the real-valued functions
\beq \label{Idef}
I_\ell(u) = \int_0^u \frac{s^\ell}{1+s^2 + 2 s \cos \alpha} \: ds\:,\qquad \ell=1,2\:.
\eeq
Choosing a polar decomposition of~$b$ and~$d$,
\[ b = |b|\, e^{i \beta}\:,\quad d = |d|\, e^{i \delta} 
\qquad \text{with~$\beta, \delta \in [-\pi, \pi]$} \:, \]
and introducing the abbreviation
\beq \label{rhodef}
\rho = \left| \frac{b}{8d} \right| ,
\eeq
we obtain
\[ \re \big( f(t)-f(t_0) \big) = |d| \:h(u) \:, \]
where
\beq \label{hdef}
h(u) := \cos(2 \alpha+\delta) \, u^2 + 2 \rho\,
\Big( \cos(2 \alpha+\beta)\, I_1(u) + \cos(\alpha+\beta)\, I_2(u) \Big) \:.
\eeq

\begin{Lemma} \label{lemma31}
The function~$h$ satisfies the inequality
\[ \left| \frac{h(u)}{u^2} - \cos(2 \alpha+\delta) +  \rho \:\frac{\sin(\alpha+\beta)}{2 \sin \alpha} \right|
\leq \frac{\rho}{2 \sin \alpha} \:. \]
\end{Lemma}
\Proof It is obvious from~\eqref{hdef} and~\eqref{Idef} that~$h(0)=0=h'(0)$, and that~$h''(0)$
exists and is finite.
Suppose that~$u>0$. Then
\[ h'(u) \leq u \sup_{s \in \R \setminus \{0\}} \frac{h'(s)}{s} \:, \]
and integrating from~$0$ to~$u$ yields
\[ \frac{h(u)}{u^2} \leq \sup_{s \in \R \setminus \{0\}} \frac{h'(s)}{2 s}\:. \]
Similarly, one can estimate~$h(u)/u^2$ from below by~$\inf h'(s)/(2s)$.
Proceeding analogously in the case~$u<0$ gives
\[ \inf_{s \in \R \setminus \{0\}} \frac{h'(s)}{2 s}
\leq \frac{h(u)}{u^2} \leq  \sup_{s \in \R \setminus \{0\}} \frac{h'(s)}{2 s}\:. \]

By differentiating~\eqref{Idef}, we readily obtain
\beq \label{hpeq}
\frac{h'(s)}{2 s} = \cos(2 \alpha+\delta) + \rho\; \frac{
\cos(2 \alpha+\beta) + s\: \cos(\alpha+\beta)}{1 + s^2 + 2 s \cos \alpha} \:,
\eeq
which implies that
\[ \lim_{s \rightarrow \pm \infty} \frac{h'(s)}{2 s} = \cos(2 \alpha+\delta)\:. \]
As a consequence, the function~$h'(s)/(2s)$ is bounded from above and below
by its inner maximum and minimum, respectively.
Differentiating~\eqref{hpeq} and computing the zeros, a straightforward computation gives the result.
\QED

\begin{Lemma} \label{lemmaphases}
By appropriately choosing~$\alpha$ in~\eqref{alphadef}, one can arrange that
\beq \label{slower}
\sin \alpha \geq \frac{1}{2}
\eeq
and
\beq \label{refes}
\re (f(t) - f(t_0)) \leq -\kappa\,|d|\: |\tau|^2 \:,
\eeq
where
\beq \label{kappaeq}
\kappa = \left\{ \begin{array}{cl} \displaystyle
1 - \rho\: \frac{\sin^2 \left(\frac{\beta}{2} - \frac{\delta}{4} \right)}{\cos \frac{\delta}{2}} & 
\text{if~$\displaystyle \delta \in \Big(\!-\frac{2 \pi}{3}, \frac{2 \pi}{3} \Big)$} \\[1em]
\displaystyle -\cos \left( \delta \pm \frac{\pi}{3} \right) - 2 \rho \left( 1 - \sin \left( \beta \pm \frac{\pi}{6} \right)
\right) & \text{if~$\displaystyle \delta \not\in \Big(\!-\frac{2 \pi}{3}, \frac{2 \pi}{3} \Big)$} \end{array} \right.
\eeq
(and one can choose the cases~$\pm$ as desired).
If~$\rho<1/8$, we can choose~$\kappa=1/4$.
\end{Lemma}
\Proof In the case~$\delta \in (-\frac{2 \pi}{3}, \frac{2 \pi}{3})$, we choose
\[ \alpha = \frac{\pi - \delta}{2}\:, \]
which gives~$\cos(2 \alpha-\delta)=-1$. Applying Lemma~\ref{lemma31} and using
the identities
\begin{align*}
\sin \left(\frac{\pi - \delta}{2} \right) &= \cos \left( \frac{\delta}{2} \right) \\
1-\sin(\alpha+\beta) &= 2 \sin^2 \left( \frac{\alpha+\beta}{2} - \frac{\pi}{4} \right)
= 2 \sin^2 \left( \frac{\beta}{2} - \frac{\delta}{4} \right)
\end{align*}
gives~\eqref{refes}.

In the second case, we apply Lemma~\ref{lemma31} choosing
\[ \alpha = \frac{\pi}{6} \qquad \text{or} \qquad \alpha = \frac{5 \pi}{6}\:. \]
A short calculation again gives~\eqref{refes}.

If~$\rho<1/8$, another short calculation  shows that~$\kappa$ as given by~\eqref{kappaeq}
is always larger than~$1/4$. Therefore, we can choose~$\kappa=1/4$.
\QED

\begin{figure}
\begin{picture}(0,0)%
\includegraphics{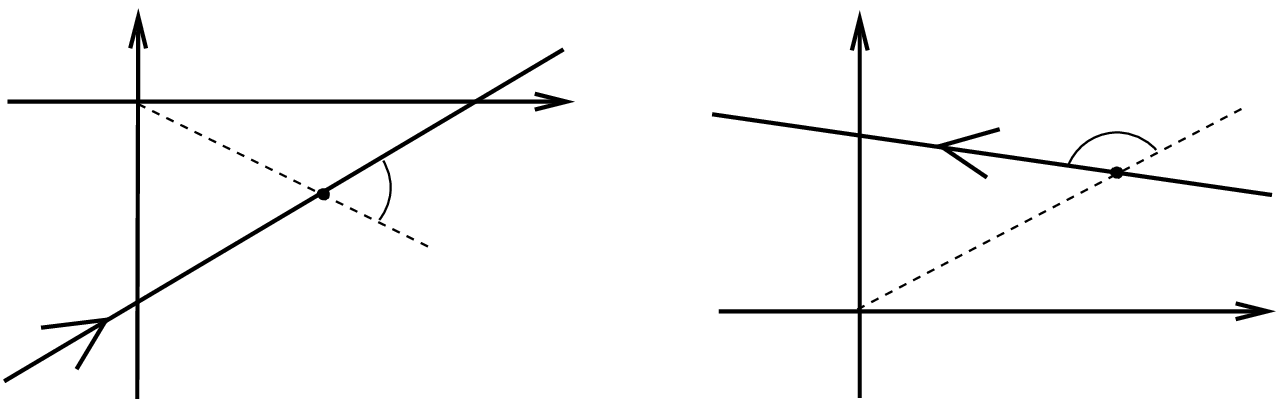}%
\end{picture}%
\setlength{\unitlength}{2486sp}%
\begingroup\makeatletter\ifx\SetFigFont\undefined%
\gdef\SetFigFont#1#2#3#4#5{%
  \reset@font\fontsize{#1}{#2pt}%
  \fontfamily{#3}\fontseries{#4}\fontshape{#5}%
  \selectfont}%
\fi\endgroup%
\begin{picture}(9755,3500)(72,-6760)
\put(3976,-4422){\makebox(0,0)[lb]{\smash{{\SetFigFont{11}{13.2}{\familydefault}{\mddefault}{\updefault}$\re t$}}}}
\put(9398,-6019){\makebox(0,0)[lb]{\smash{{\SetFigFont{11}{13.2}{\familydefault}{\mddefault}{\updefault}$\re t$}}}}
\put(2364,-4541){\makebox(0,0)[lb]{\smash{{\SetFigFont{11}{13.2}{\familydefault}{\mddefault}{\updefault}$t_0$}}}}
\put(3113,-4767){\makebox(0,0)[lb]{\smash{{\SetFigFont{11}{13.2}{\familydefault}{\mddefault}{\updefault}$\alpha$}}}}
\put(1313,-3559){\makebox(0,0)[lb]{\smash{{\SetFigFont{11}{13.2}{\familydefault}{\mddefault}{\updefault}$\im t$}}}}
\put(947,-6664){\makebox(0,0)[lb]{\smash{{\SetFigFont{11}{13.2}{\familydefault}{\mddefault}{\updefault}$\arg(e^{i \alpha}\, t_0) \in (-\frac{\pi}{4}, \frac{\pi}{4})$}}}}
\put(6750,-3530){\makebox(0,0)[lb]{\smash{{\SetFigFont{11}{13.2}{\familydefault}{\mddefault}{\updefault}$\im t$}}}}
\put(8362,-4152){\makebox(0,0)[lb]{\smash{{\SetFigFont{11}{13.2}{\familydefault}{\mddefault}{\updefault}$\alpha$}}}}
\put(6451,-6650){\makebox(0,0)[lb]{\smash{{\SetFigFont{11}{13.2}{\familydefault}{\mddefault}{\updefault}$\arg(e^{i \alpha} \,t_0) \in (\frac{3 \pi}{4}, \frac{5 \pi}{4})$}}}}
\put(1622,-5684){\makebox(0,0)[lb]{\smash{{\SetFigFont{11}{13.2}{\familydefault}{\mddefault}{\updefault}$t(u) = t_0\, (1+e^{i \alpha} u)$}}}}
\end{picture}%
\caption{Contour deformations of~$\Gamma_+$ and~$\Gamma_-$.}
\label{figdeform}
\end{figure}

\begin{Thm} \label{thm1} Suppose that the phases~$\beta$ and~$\delta$ are such that the
parameter~$\kappa$ in Lemma~\ref{lemmaphases} is positive. We consider the
parabolic cylinder function
\[ U_a(z) = \left\{ \begin{array}{cl} U_a^+(z) & \text{if~$\re(t_0 e^{i \alpha}) < 0$} \\[0.2em]
U_a^-(z) & \text{if~$\re(t_0 e^{i \alpha}) \geq 0$} \:. \end{array} \right. \]
Then this parabolic cylinder function can be approximated by the function
\beq \label{Utilde}
\tilde{U}_a(z)  = e^{-\frac{z^2}{4} + f(t_0)} \:t_0\: \frac{2 \sqrt{2 \pi}}{\sqrt{-8d-b}} \:,
\eeq
with the relative error bounded by
\[ \left| \frac{U_a(z)}{\tilde{U}_a(z)} - 1 \right| 
\leq \frac{\sqrt{|8d+b|}\, |b|}{\kappa^2\:|d|^2} \left( 1
+ \log^2 \bigg( 1+\frac{(\kappa \, |d|)^\frac{3}{2}}{|b|} \bigg) \right) . \]
\end{Thm}
\Proof We first consider the contour deformations in more detail.
It is guaranteed by Lemma~\ref{lemmaphases} that~$\re(t^2)$ is positive
at both ends of the contour. In view of~\eqref{tpar} and~\eqref{alphadef}, this implies that
\[ \arg \big( e^{i \alpha} \, t_0 \big) \in \Big( -\frac{\pi}{4}, \frac{\pi}{4} \Big) \cup
\Big( \frac{3\pi}{4}, \frac{5\pi}{4} \Big) \mod 2 \pi \:. \]
In the case~$\arg ( e^{i \alpha}\,  t_0 ) \in ( -\frac{\pi}{4}, \frac{\pi}{4} )$,
the contour can be continuously deformed to~$\Gamma_-$ (see Figure~\ref{figdeform}).
Likewise, in the case~$\arg ( e^{i \alpha}  t_0 ) \in ( \frac{3 \pi}{4}, \frac{5 \pi}{4})$,
a contour deformation gives the contour~$\Gamma_+$, but with the opposite orientation
(see again Figure~\ref{figdeform}).
Hence from~\eqref{Uaz}, \eqref{tpar} and~\eqref{alphadef}, we obtain
\beq \label{Uax}
U_a(z) = \pm e^{-\frac{z^2}{4} + i \alpha + f(t_0)} \:t_0 \int_{-\infty}^\infty e^{f(t)-f(t_0)}\: du\:.
\eeq
Replacing~$f$ by~$\tilde{f}$, we obtain a Gaussian integral, which can easily be computed
to obtain~\eqref{Utilde}, with an appropriate choice of the sign of the square root~$\sqrt{-8d-b}$.

In order to estimate the integral for large~$|u|$, we fix a parameter~$L>0$. Then
\begin{align*}
\bigg| \int_{\R \setminus [-L,L]} &e^{f(t)-f(t_0)}\: du \bigg|
= \int_{\R \setminus [-L,L]}e^{\re(f(t)-f(t_0))}\: du 
\overset{\eqref{refes}}{\leq}
\int_{\R \setminus [-L,L]} e^{-\kappa\, |d|\, u^2}\: du \\
&\leq e^{-\frac{1}{2}\, \kappa\, |d|\, u^2}
\int_{-\infty}^\infty e^{-\frac{1}{2}\, \kappa\, |d|\, u^2}\: du
= \sqrt{\frac{2 \pi}{\kappa\, |d|}}\: e^{-\frac{1}{2} \,\kappa\, |d|\, L^2}\:.
\end{align*}
Since~$\tilde{f}$ is the quadratic approximation to~$f$ at~$\tau=0$,
the inequality~\eqref{refes} also holds for~$\tilde{f}$.
Hence the last estimate is also true for~$\tilde{f}$,
\[ \bigg| \int_{\R \setminus [-L,L]} e^{\tilde{f}(t)-f(t_0)}\: du \bigg|
\leq \sqrt{\frac{2 \pi}{\kappa\, |d|}}\: e^{-\frac{1}{2} \,\kappa\, |d|\, L^2}\:. \]

Next, on the interval~$[-L,L]$ we use the estimates
\beql{smalles}
\begin{split}
\Big| e^{f(t)-f(t_0)} &- e^{\tilde{f}(t)-f(t_0)} \Big| = 
\bigg| \int_0^1 \frac{d}{ds} e^{s f(t) + (1-s) \tilde{f}(t)-f(t_0)}\: ds \bigg| \\
&= \int_0^1 \left| (f(t)-\tilde{f}(t)) \:e^{s f(t) + (1-s) \tilde{f}(t)-f(t_0)} \right| ds \\
&\leq |f(t) - \tilde{f}(t)|\: \sup_{s \in [0,1]} e^{\re \big( s f(t) + (1-s) \tilde{f}(t)-f(t_0) \big)} 
\leq |f(t) - \tilde{f}(t)|\:,
\end{split}
\eeq
where in the last step we used that the real parts of~$f(t)-f(t_0)$
and~$\tilde{f}(t) - \tilde{f}(t_0)$ are both negative in view of~\eqref{refes}
and the fact that~$\tilde{f}$ is the quadratic Taylor polynomial of~$f$ about~$t_0$.

Moreover, as~$\tilde{f}$ is the quadratic approximation to~$f$, we know that
\beq
|f(t) - \tilde{f}(t)| \leq \frac{|\tau|^3}{3!}\: \sup_{\tau} |f^{(3)}(\tau)|
\overset{\eqref{ft0rel}}{\leq} \frac{|\tau|^3}{3!}\: \frac{|b|}{4} \sup_{\tau} \frac{2}{|1+\tau|^3} \:.
\label{quadratic}
\eeq
Since the distance from our contour~\eqref{alphadef} to the point~$-1$
is equal to~$\sin \alpha$, we can apply~\eqref{slower} to obtain
\[ |f(t) - \tilde{f}(t)| \leq \frac{2}{3}\: |b|\, |\tau|^3\:. \]
We conclude that
\[ \bigg| \int_{-L}^L \left( e^{f(t)-f(t_0)} - e^{\tilde{f}(t)-f(t_0)} \right) du \bigg|
\leq \frac{2}{3}\: |b| \int_{-L}^L |u|^3\: du = \frac{1}{3}\: |b| \,L^4 \]

Combining the above estimates, we conclude that for any~$L>0$, the following inequality holds,
\beq \label{UUt}
\bigg| \int_{-\infty}^\infty \left( e^{f(t)-f(t_0)} - e^{\tilde{f}(t)-f(t_0)} \right) du \bigg| \leq
2 \:\sqrt{\frac{2 \pi}{\kappa\, |d|}}\: e^{-\frac{1}{2} \,\kappa\, |d|\, L^2}
+ \frac{1}{3}\: |b|\, L^4 \:.
\eeq
We want to choose~$L$ in such a way that the two error terms are of comparable size.
To this end, we set
\beq \label{Lchoice2}
L = \sqrt{\frac{2}{\kappa\, |d|}}\, \log^\frac{1}{2} \bigg( 1 + \frac{(\kappa \, |d|)^\frac{3}{2}}{|b|} \bigg)
\eeq
to obtain
\begin{align}
\bigg| \int_{-\infty}^\infty &\left( e^{f(t)-f(t_0)} - e^{\tilde{f}(t)-f(t_0)} \right) du \bigg| \nonumber \\
&\leq 2 \:\sqrt{\frac{2 \pi}{\kappa\, |d|}}\: \frac{|b|}{(\kappa \, |d|)^\frac{3}{2}}
+ \frac{|b|}{3}\: \frac{4}{\kappa^2\, |d|^2}  \log^2 \bigg( 1+\frac{(\kappa \, |d|)^\frac{3}{2}}{|b|} \bigg) \nonumber \\
&= \frac{|b|}{\kappa^2\, |d|^2} \left( 2\: \sqrt{2 \pi}
+ \frac{4}{3} \log^2 \bigg( 1+\frac{(\kappa \, |d|)^\frac{3}{2}}{|b|} \bigg) \right) . \label{intes}
\end{align}
Computing the numerical constants and using~\eqref{Uax} and~\eqref{Utilde} gives the result.
\QED

\begin{Thm} \label{thm2} Under the assumptions of Theorem~\ref{thm1},
\begin{align*}
\left| \frac{d}{dz} \bigg( \frac{U_a(z)}{\tilde{U}_a(z)} \bigg) \right| &\leq 
\frac{|b| \,|d'(z)|}{\kappa^2 |d|^2 \sqrt{|8d+b|}} \left(1 + \log^2 \bigg( 1+\frac{(\kappa \, |d|)^\frac{3}{2}}{|b|}
\bigg) \right) \\
&\quad + \frac{|b| |d'(z)|}{\kappa^3 |d|^3} \:\sqrt{|8d+b|} \left(1 + \log^3 \bigg( 1+\frac{(\kappa \, |d|)^\frac{3}{2}}
{|b|} \bigg) \right).
\end{align*}
\end{Thm}
\Proof Introducing the abbreviations
\[ I(z) = \int_{-\infty}^\infty e^{f(t)-f(t_0)}\: du \:,\qquad
\tilde{I}(z) = \int_{-\infty}^\infty e^{\tilde{f}(t)-f(t_0)}\: du\:, \]
we have
\[ \frac{d}{dz} \left( \frac{U_a(z)}{\tilde{U}_a(z)} \right)
= \frac{d}{dz} \left( \frac{I(z)}{\tilde{I}(z)} \right)
= \frac{I'(z) - \tilde{I}'(z)}{\tilde{I}(z)} + \frac{\tilde{I}'(z)}{\tilde{I}(z)}\:
\frac{\tilde{I}(z) - I(z)}{\tilde{I}(z)}\:. \]
Computing the Gaussian integral gives
\[ \tilde{I}(z) = \left( \frac{8 \pi}{-8d-b} \right)^{\frac{1}{2}}\:, \]
and a straightforward calculation using~\eqref{intes} yields
\[ \left| \frac{\tilde{I}'(z)}{\tilde{I}(z)}\:
\frac{\tilde{I}(z) - I(z)}{\tilde{I}(z)} \right|
\leq \frac{\sqrt{2}}{3 \:\sqrt{\pi}}\: \frac{|b|}{\kappa^2 |d|^2\, \sqrt{|8d+b|}}\:
|d'(z)| \left( 1 + \log^2 \bigg( 1+\frac{(\kappa \, |d|)^\frac{3}{2}}{|b|} \bigg) \right) . \]

Moreover,
\[ I'(z)-\tilde{I}'(z) = \int_{-\infty}^\infty d'(z) \:\tau(u)^2 \left( e^{f(t)-f(t_0)} - e^{\tilde{f}(t)-f(t_0)} \right) du \]
and thus
\[ |I'(z) - \tilde{I}'(z)| \leq |d'(z)|
\bigg| \int_{-\infty}^\infty u^2 \left( e^{f(t)-f(t_0)} - e^{\tilde{f}(t)-f(t_0)} \right) du \bigg| \]
Proceeding as in the proof of Theorem~\ref{thm1}, we obtain similar to~\eqref{UUt}
\begin{align*}
\bigg| \int_{-\infty}^\infty u^2 \left( e^{f(t)-f(t_0)} - e^{\tilde{f}(t)-f(t_0)} \right) du \bigg|
\leq \frac{2 \:\sqrt{2 \pi}}{(\kappa\, |d|)^\frac{3}{2}}\: e^{-\frac{1}{2} \,\kappa\, |d|\, L^2}
+ \frac{2}{9}\: |b|\, L^6 \:.
\end{align*}
Again choosing~$L$ according to~\eqref{Lchoice2}, we obtain
\[ \bigg| \int_{-\infty}^\infty u^2 \left( e^{f(t)-f(t_0)} - e^{\tilde{f}(t)-f(t_0)} \right) du \bigg|
\leq \frac{|b|}{\kappa^3 \,|d|^3} \left( 2 \sqrt{2 \pi} + \frac{16}{9}\:
\log^3 \bigg( 1+\frac{(\kappa \, |d|)^\frac{3}{2}}{|b|} \bigg) \right) . \]
Computing the numerical constants and combining all the terms gives the result. 
\QED

We finally explain how the parabolic cylinder function is related to
the WKB wave function~\eqref{parWKB}.
\begin{Remark} {\em{ Assume that~$|z|^2 > 4 \,|b|$.
Choosing~$t_0=t_+$, we know from~\eqref{cdef2} that
\begin{align*}
\left| \sqrt{z^2-b} - z \right| &\leq |b| \sup_{4 |\zeta| < |z|^2} \left| \partial_\zeta \big(\sqrt{z^2 - \zeta} \big) \right|
\leq |b| \sup_{4 |\zeta| < |z|^2} \;\frac{1}{2 \left| \sqrt{z^2 - \zeta} \right|}
\leq \frac{|b|}{|z|} \\
|d| &\geq \frac{1}{8} \left( |2z| - \left| \sqrt{z^2-b} - z \right| \right)^2 
\geq \frac{|z|^2}{2} \left( 1 - \frac{|b|}{2 |z|^2} \right)^2
> \frac{|z|^2}{4} \\
\rho \;&\!\!\!\overset{\eqref{rhodef}}{<}
\frac{|b|}{2 |z|^2} < \frac{1}{8} \:.
\end{align*}
Hence choosing~$\kappa=1/4$, Theorems~\ref{thm1} and~\ref{thm2} apply,
showing that~$U_a(z)$ is well-approximated by~$\tilde{U}_a(z)$ with a small error.
Next, a straightforward calculation yields
\beq \label{Uatil}
\tilde{U}_a = 2^\frac{b}{4} \,e^{\frac{b}{8} + i \alpha}\: \sqrt{\pi}\;
\frac{\big( z + \sqrt{z^2-b}\, \big)^{\frac{1}{2}-\frac{b}{4}}}{(z^2-b)^\frac{1}{4}}\:
\exp \left( {\frac{z}{4}\: \sqrt{z^2-b}} \right) \:.
\eeq
Thus we see that replacing~$z^2-b$ by~$z^2-b+2$, we obtain up to a constant,
the WKB wave function~\eqref{parWKB}.
Assuming furthermore that~$|z|$ is sufficiently large, we conclude that~$U_a(z)$
is indeed well-approximated by~$\phi_\WKB(z)$.
}} \QEDrem
\end{Remark}

In particular, this last remark together with the identity~\eqref{Uatil} enables us to compute the following limit.
\begin{Corollary} \label{cor36} For any~$z_0 \in \C$,
\[ \lim_{t \rightarrow \infty} \frac{1}{t}\: \frac{U_a'(t z_0)}{U_a(t z_0)} = \frac{z_0}{2} \:. \]
\end{Corollary}

\subsection{Ruling out Zeros of the Parabolic Cylinder Functions} \label{sec32}
In analogy to Section~\ref{secairyzero}, we now analyze the solution~$U_a(z)$ along the
straight line~\eqref{zudef}. Thus setting~$\phi(u) = U_a(z(u))$, this function
is a solution of the Sturm-Liouville equation~\eqref{5ode} with
\[ V(u) = \lambda^2 \left( \frac{(z_0 + \lambda u)^2}{4}+a \right) . \]
After suitably choosing~$\lambda$, one can apply Corollary~\ref{cor25} to obtain the following result.
\begin{Thm} The function~$U_a$ has no zeros in the complex plane.
\end{Thm}
\Proof We always choose~$\lambda$ such that~$\lambda^4$ is real. We then obtain
\[ \im V(u) = \frac{1}{4}\,
\im \Big( 2 \lambda^3 z_0\, u + \lambda^2 \left( z_0^2 + 4 a \right) \Big) . \]
Moreover, according to Corollary~\ref{cor36}, we know that the function~$y(u) := \phi'(u)/\phi(u)$
behaves asymptotically for large negative~$u$ as
\[ y(u) = \lambda\: \frac{U_a'(z(u))}{U_a(z(u))} \sim
\lambda\: \frac{z_0 + \lambda u}{2} \:. \]

We first consider the case~$\re z_0 \neq 0$. Choosing~$\lambda=\pm i$, we get
\begin{align*}
\im y(u) &\sim \pm \frac{1}{2} \,\re z_0 \\
\im V(u) &= \frac{1}{4}\:
\im \Big(\mp 2 i z_0\, u - \left( z_0^2 + 4 a \right) \Big)
= \mp \frac{u}{2} \,\re z_0 - \im \left( \frac{z_0^2}{4} + a \right) .
\end{align*}
Thus in the limit~$u \rightarrow -\infty$, the imaginary parts of the functions~$y$ and~$V$
are non-zero and have the same sign. Moreover, by choosing the sign~$\pm$ appropriately,
we can arrange that the imaginary part of~$V$ does not change sign on the interval~$(-\infty, 0]$.
Exactly as in the proof of Proposition~\ref{thmairy2} it follows that~$U_a(z_0) \neq 0$.

In the remaining case~$\re z_0 = 0$, we choose~$\lambda = \sigma (1 + i \tau)$
with~$\sigma, \tau \in \{\pm 1\}$. Then~$\lambda^2 = 2 i \tau$ and~$\lambda^4$ is indeed real.
It follows that
\begin{align*}
\im y(u) &\sim \frac{u}{2} \: \im ( \lambda^2 ) = \tau u \\
4 \im V(u) &=
\im \Big( 2 \lambda^3 z_0\, u + \lambda^2 \left( z_0^2 + 4 a \right) \Big) \\
&= 2 \tau \re \Big( 2 \lambda z_0\, u + \left( z_0^2 + 4 a \right) \Big) \\
&= -4 \sigma u \im (z_0) + 2 \tau \re  \big( z_0^2 + 4 a \big) \:.
\end{align*}
Thus in the limit~$u \rightarrow -\infty$, the imaginary part of~$y$ is non-zero, and
by choosing~$\sigma$ appropriately we can arrange that it has the same sign as
the imaginary part of~$V$. Moreover, by adjusting~$\tau$ we can arrange that~$\im V$
does not change sign on the interval~$(-\infty, 0]$. Again proceeding exactly as
in the proof of Proposition~\ref{thmairy2}, we obtain that~$U_a(z_0) \neq 0$.
\QED

\subsection{The Airy Limit of Parabolic Cylinder Functions}  \label{secparairy}
In this section, we identify the parameter range where the parabolic cylinder function
is well-approximated by the Airy function.
The idea is to expand the function~$f$ in~\eqref{Uaz} around a point~$t_0$ where
the second derivative vanishes. Then the Taylor polynomial of~$f$ around~$t_0$
involves a linear and a cubic term, just like the exponent in the integral representation~\eqref{Arep}
of the Airy function. The zeros of~$f''(t)$ are computed to be~$\pm \sqrt{b}/2$.
In order to fix the sign, we let~$c$ be a solution of the equation
\[ c^6 = \frac{b}{4} \]
and set
\[ t_0 = c^3\:. \]
Then the function~$f$ has the Taylor expansion
\[ f(t) = f(t_0) + \big( z - 2 c^3 \big) (t-t_0) - \frac{1}{3 c^3} \:(t-t_0)^3 + \O \big( (t-t_0)^4 \big) . \]
Defining~$\tilde{f}$ as the corresponding Taylor polynomial and introducing the parametrization
\beq \label{taupar}
t = t_0 - c \tau \:,
\eeq
we obtain
\begin{align*}
f(t) - f(t_0) &= \left( c^4 - c z \right) \tau - \frac{c^2}{2}\: \tau^2 - c^6 \:  \log \Big(1- \frac{\tau}{c^2} \Big) \\
\tilde{f}(t) - f(t_0) &= \left( 2 c^4 - c z \right) \tau + \frac{\tau^3}{3} \:.
\end{align*}
Then writing~\eqref{Uaz} as
\beq \label{Uaair}
U_a(z) = e^{-\frac{z^2}{4}} \int_\Gamma e^{\tilde{f}(t)}\: dt 
+ e^{-\frac{z^2}{4}} \int_\Gamma \big( e^{\tilde{f}(t)}- e^{f(t)} \big)\: dt  \:,
\eeq
the first integral is just of the form of the Airy function~\eqref{Arep},
whereas the second integral is the error term.
In order to show that the parabolic cylinder function goes over to the
Airy function, we must make sure that the contour~$\Gamma$ is compatible
with that in~\eqref{Arep}, and that the functions~$e^{\tilde{f}(t)}$ and~$e^{f(t)}$
both decay exponentially on both ends of the contour.

\begin{figure}
\begin{picture}(0,0)%
\includegraphics{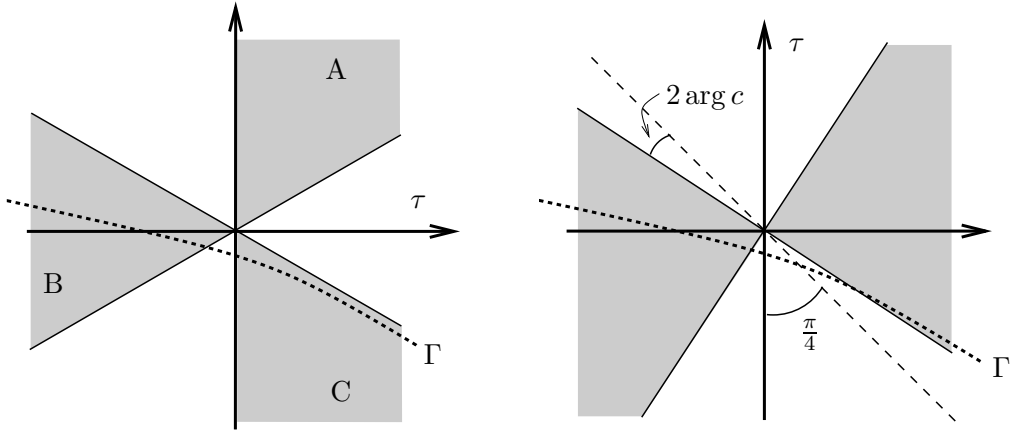}%
\end{picture}%
\setlength{\unitlength}{2486sp}%
\begingroup\makeatletter\ifx\SetFigFont\undefined%
\gdef\SetFigFont#1#2#3#4#5{%
  \reset@font\fontsize{#1}{#2pt}%
  \fontfamily{#3}\fontseries{#4}\fontshape{#5}%
  \selectfont}%
\fi\endgroup%
\begin{picture}(9892,4342)(-110,-7470)
\put(3099,-3942){\makebox(0,0)[lb]{\smash{{\SetFigFont{11}{13.2}{\familydefault}{\mddefault}{\updefault}A}}}}
\put(3151,-7152){\makebox(0,0)[lb]{\smash{{\SetFigFont{11}{13.2}{\familydefault}{\mddefault}{\updefault}C}}}}
\put(7719,-3664){\makebox(0,0)[lb]{\smash{{\SetFigFont{11}{13.2}{\familydefault}{\mddefault}{\updefault}$\tau$}}}}
\put(6496,-4196){\makebox(0,0)[lb]{\smash{{\SetFigFont{11}{13.2}{\familydefault}{\mddefault}{\updefault}$2 \arg c$}}}}
\put(3946,-5232){\makebox(0,0)[lb]{\smash{{\SetFigFont{11}{13.2}{\familydefault}{\mddefault}{\updefault}$\tau$}}}}
\put(287,-6080){\makebox(0,0)[lb]{\smash{{\SetFigFont{11}{13.2}{\familydefault}{\mddefault}{\updefault}B}}}}
\put(9751,-6912){\makebox(0,0)[lb]{\smash{{\SetFigFont{11}{13.2}{\familydefault}{\mddefault}{\updefault}$\Gamma$}}}}
\put(4081,-6784){\makebox(0,0)[lb]{\smash{{\SetFigFont{11}{13.2}{\familydefault}{\mddefault}{\updefault}$\Gamma$}}}}
\put(7830,-6544){\makebox(0,0)[lb]{\smash{{\SetFigFont{11}{13.2}{\familydefault}{\mddefault}{\updefault}$\frac{\pi}{4}$}}}}
\end{picture}%
\caption{Admissible regions for the integration contour.}
\label{figregions}
\end{figure}

We now demonstrate how to continuously deform the original integration contour~$t \in \R + i$ in such a way
that the real parts of both~$f$ and~$\tilde{f}$ are both negative at both ends of the new contour.
In our parametrization~\eqref{taupar}, we need to specify the contour in the variable~$\tau$.
As only the product~$c \tau$ enters in~\eqref{taupar}, we can assume without loss of generality that
\[ \frac{3 \pi}{4} \leq \arg c < \frac{13 \pi}{12} \mod 2 \pi\:. \]
Then we can deform the contour such that in Figure~\ref{figregions}
the integration begins in region~B and ends in region~C
(we could also end in region~A, but our choice fits well with the contour
chosen for the Airy function in~\eqref{Arep}).
For the real part of~$-c^2 \tau^2$ to be negative at both ends of the contour, we
need to assume that (see Figure~\ref{figregions})
\[ -\frac{5 \pi}{12} < 2 \arg c < \frac{\pi}{12} \mod 2 \pi \:. \]
Provided that this inequality holds, there is indeed a contour such that the
real parts of both~$f$ and~$\tilde{f}$ are negative at both ends of the contour
(see Figure~\ref{figregions}). Then the error term in~\eqref{Uaair} decays exponentially
on both ends of the contour.

Having accomplished this exponential decay, rigorous error estimates
could be obtained using methods similar to those
in Sections~\ref{secairyWKB} and~\ref{secparWKB}.
Giving the details would go beyond the scope of this paper.

\subsection{The Airy-WKB Limit of Parabolic Cylinder Functions}
We now consider an asymptotics which includes both the WKB and Airy asymptotics
of Sections~\ref{secparWKB} and~\ref{secparairy}. Our main result is to obtain
approximate solutions in terms of Airy functions and to derive rigorous error bounds.
As in Section~\ref{secparWKB}, we expand the function~$f(t)$ in the exponent of the
integral representation~\eqref{Uaz} in a Taylor series around~$t_0 = t_\pm$.
But now we choose~$\tilde{f}$ as the cubic Taylor approximation.
Thus according to~\eqref{f23} and~\eqref{tpar}, we set
in analogy to~\eqref{ft0rel} and~\eqref{ftt0rel},
\begin{align}
f(t) - f(t_0) &= -\frac{t_0^2}{2}\: \tau^2 + \frac{b}{4} \,\big( \tau - \log(1+\tau) \big) \label{f} \\
\tilde{f}(t) - f(t_0) &= \left(-\frac{t_0^2}{2} + \frac{b}{8} \right) \tau^2 
- \frac{b}{12}\: \tau^3\:. \label{ftt}
\end{align}
We introduce the angle~$\gamma$ by
\[ -\frac{b}{t_0^3} = \left| \frac{b}{t_0^3} \right| e^{3 i \gamma}
\qquad \text{with} \qquad
\gamma \in \big(-60^\circ, 60^\circ]\:. \]
Using an argument similar to that shown in Figure~\ref{figregions}, one can verify that if
\beql{gamrange}
\gamma \in (-60^\circ, 0) \:,
\eeq
then the contour~$\Gamma$ can be deformed in such a way that the following conditions hold:
\begin{itemize}
\item[(a)] The real parts of both~\eqref{f} and~\eqref{ftt} decay at both ends of the contour.
\item[(b)] The contour is (up to a rotation) a deformation of the Airy contour in~\eqref{Arep}.
\end{itemize}
We then introduce the approximate solution~$\tilde{U}_a$ by
\[ \tilde{U}_a(z) = e^{-\frac{z^2}{4}} \int_{\Gamma} e^{\tilde{f}(t)}\: dt \:. \]
As the exponent is a cubic polynomial~\eqref{ftt}, the approximate solution can be
expressed explicitly in terms of the Airy function.
\begin{Lemma} For the above choice of the contour,
\beq \label{UaAiry}
\tilde{U}_a(z) = 2^\frac{2}{3}\: \sqrt{\pi}\:
|t_0|\: |b|^{-\frac{1}{3}} \: t_0^{-\frac{b}{4}}\:
e^{-\frac{2}{3}\: h(z)^3 -\frac{t_0^2}{2} + t_0 z - \frac{z^2}{4} - i \gamma}\: \A \big( h(z)^2 \big)\:,
\eeq
where
\beq \label{hdef2}
h(z) := \frac{e^{-2 i \gamma}}{2 \cdot 2^{2/3}}\: \frac{|t_0|^2}{t_0^2}\:
|b|^{-\frac{2}{3}}\:(4 t_0^2-b)  \:.
\eeq
\end{Lemma}
\Proof We parametrize the deformed contour by
\[ \tau(u) = \left\{ \begin{array}{cl} e^{i (-\arg t_0 + \beta_+)}\,u & \text{if~$u \geq 0$} \\
e^{i (-\arg t_0 + \beta_-)}\,u & \text{if~$u < 0$\:.} \end{array} \right. \]
Then
\[ \tilde{f}(t) - f(t_0) = \left(-\frac{|t_0|^2}{2} + \frac{b}{8}\, \frac{|t_0|^2}{t_0^2} \right) e^{2 i \beta_\pm} \,u^2
+ \frac{|b|}{12} \left( e^{i (\gamma + \beta_\pm)} \:u \right)^3\:. \]
In order to compute the corresponding contour integral, we shift the integration variable 
according to
\[ u \rightarrow u - \frac{|t_0|^3}{t_0^3} \:\frac{b-4 t_0^2}{2\,|b|}\: e^{-i \beta_\pm-3 i \gamma} \]
by a complex number. Then the quadratic term vanishes, so that we get an integral
of the Airy form~\eqref{Arep}. A straightforward calculation gives the result. % siehe [pawkb.nb]
\QED

By applying Theorem~\ref{thmairy} to~\eqref{UaAiry}, one can recover the WKB asymptotics,
as we now explain. Assume that~$|z| \ll |b|$.
According to~\eqref{tpm}, we can choose~$t_0=t_\pm$ such that~$|t_0| \ll |b|$.
Then the argument of~$h(z)$ is given asymptotically by
\[ \lim_{|t_0| \rightarrow \infty} \arg h(z) \overset{\eqref{hdef2}}{=}
\lim_{|t_0| \rightarrow \infty} \arg
\frac{e^{-2 i \gamma}}{2 \cdot 2^{2/3}}\: 4\, |b|^{-\frac{2}{3}}
= -2 \gamma \overset{\eqref{gamrange}}{\in} (0, 120^\circ) \:. \]
Again defining the roots by~$z^\alpha = \exp(\alpha \log(z))$
with the branch cut along the ray~\eqref{cut}, one sees that
for large enough~$|t_0|$, the root of~$h^2$ is given by~$\sqrt{h^2} =+h$.
Using this fact in the WKB wave function~\eqref{AWKB}, we obtain the asymptotics
\begin{align*}
\tilde{U}_a &\approx
2^\frac{2}{3} \:\sqrt{\pi}\: \frac{|t_0|}{\sqrt{h}\: |b|^\frac{1}{3}}\: \exp \left(
-\frac{t_0^2}{2}+t_0 \,z -\frac{z^2}{4} -\frac{b}{4} -i \gamma \right) \\
&= 2\:\sqrt{2 \pi}\: \frac{t_0}{\sqrt{4 t_0^2 - b}} \exp \left( -\frac{t_0^2}{2}+t_0 \,z -\frac{z^2}{4} -\frac{b}{4}\: \log t_0 \right) ,
\end{align*}
where in the last step we used~\eqref{hdef2}. Keeping in mind that
\[ f(t_0) = -\frac{t^2}{2} + t_0\, z -\frac{b}{4}\: \log t_0 \qquad \text{and} \qquad
t_0^2 \overset{\eqref{cdef2}}{=} -2 d\:, \]
we get complete agreement with the WKB-approximation~\eqref{Utilde}.

The remaining task is to estimate the error of the approximation~\eqref{UaAiry}.
For simplicity, we work out these estimates only in the specific case needed
for our applications in~\cite{tinvariant, tspectral}, although our methods extend in a straightforward
way to more general cases. 
\begin{Thm} Assume that the parameters~$b$ and~$t_0$ satisfy the conditions
\begin{align}
\arg b &\in (88^\circ, 92^\circ) \label{barg} \\
\frac{4 t_0^2}{b} &\in {\mathfrak{R}} \:, \label{Rreg}
\end{align}
where the wedge-shaped region~${\mathfrak{R}} \subset \C$ is defined by (see Figure~\ref{figR})
\begin{figure}
\begin{picture}(0,0)%
\includegraphics{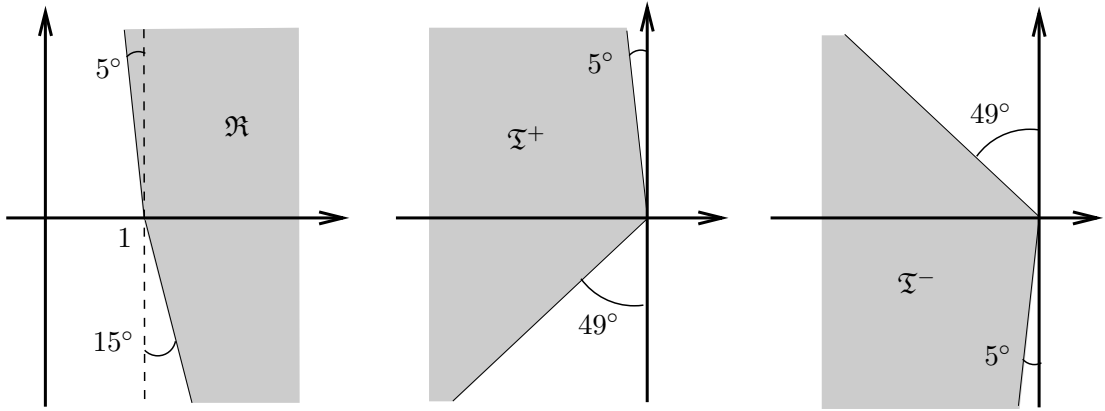}%
\end{picture}%
\setlength{\unitlength}{2486sp}%
\begingroup\makeatletter\ifx\SetFigFont\undefined%
\gdef\SetFigFont#1#2#3#4#5{%
  \reset@font\fontsize{#1}{#2pt}%
  \fontfamily{#3}\fontseries{#4}\fontshape{#5}%
  \selectfont}%
\fi\endgroup%
\begin{picture}(11006,4123)(-2792,-7384)
\put(-1649,-5749){\makebox(0,0)[lb]{\smash{{\SetFigFont{11}{13.2}{\familydefault}{\mddefault}{\updefault}$1$}}}}
\put(-1867,-4068){\makebox(0,0)[lb]{\smash{{\SetFigFont{11}{13.2}{\familydefault}{\mddefault}{\updefault}$5^\circ$}}}}
\put(-1896,-6791){\makebox(0,0)[lb]{\smash{{\SetFigFont{11}{13.2}{\familydefault}{\mddefault}{\updefault}$15^\circ$}}}}
\put(-577,-4647){\makebox(0,0)[lb]{\smash{{\SetFigFont{11}{13.2}{\familydefault}{\mddefault}{\updefault}${\mathfrak{R}}$}}}}
\put(2243,-4775){\makebox(0,0)[lb]{\smash{{\SetFigFont{11}{13.2}{\familydefault}{\mddefault}{\updefault}${\mathfrak{T}}^+$}}}}
\put(2948,-6617){\makebox(0,0)[lb]{\smash{{\SetFigFont{11}{13.2}{\familydefault}{\mddefault}{\updefault}$49^\circ$}}}}
\put(6856,-4517){\makebox(0,0)[lb]{\smash{{\SetFigFont{11}{13.2}{\familydefault}{\mddefault}{\updefault}$49^\circ$}}}}
\put(6136,-6215){\makebox(0,0)[lb]{\smash{{\SetFigFont{11}{13.2}{\familydefault}{\mddefault}{\updefault}${\mathfrak{T}}^-$}}}}
\put(3038,-4015){\makebox(0,0)[lb]{\smash{{\SetFigFont{11}{13.2}{\familydefault}{\mddefault}{\updefault}$5^\circ$}}}}
\put(6998,-6932){\makebox(0,0)[lb]{\smash{{\SetFigFont{11}{13.2}{\familydefault}{\mddefault}{\updefault}$5^\circ$}}}}
\end{picture}%
\caption{Regions in the complex plane.}
\label{figR}
\end{figure}
\[ {\mathfrak{R}} = \big\{ x+i y \:|\: x > \max \big( 1 - y \sin 5^\circ, 1 - y \sin 15^\circ \big) \big\} \:. \]
Then the parabolic cylinder function~$U_a^-$ is approximated by the
Airy function~\eqref{UaAiry}, where the following two error estimates hold,
\begin{align}
\left| U_a^-(z) - \tilde{U}_a(z) \right| 
&\leq 40 \:|t_0| \:e^{-\frac{\re z^2}{4} + f(t_0)}
\:|b|^{-\frac{2}{3}} \left( 1 + 100\, \log^\frac{5}{3} \big( 1+|b|^\frac{1}{3} \big) \right) \label{eb1} \\
\left| U_a^-(z) - \tilde{U}_a(z) \right| &\leq
20 \,|t_0| \:e^{-\frac{\re z^2}{4} + f(t_0)} \nonumber \\
&\qquad \times
\frac{|b|}{|8d+b|^\frac{5}{2}} \left( 1 + 10^5 \log^\frac{5}{2} \Big( 1 + \frac{|8d+b|^2}{|b|} \Big) \right) .
\label{eb2}
\end{align}
\end{Thm}
The remainder of this section is devoted to the proof of this theorem.
We choose
\beql{alphaval}
\alpha_- = -50^\circ = -\frac{5 \pi}{18} \qquad \text{and} \qquad \alpha_+ = 108^\circ = \frac{3 \pi}{5}\:.
\eeq
Then the contour~$\Gamma_-$ introduced after~\eqref{Uapm} can be
continuously deformed to the contour
\[ \tau \in e^{i \alpha_-} \R^+ \cup e^{i \alpha_+} \R^+ \:. \]
Moreover, this contour is a continuous deformation of the Airy contour in~\eqref{Arep}.
Hence
\beql{Umint}
U^-_a(z) = e^{-\frac{z^2}{4}} \:t_0 \left( \int_0^\infty e^{f(t)}\: e^{i \alpha_+} du
- \int_0^\infty e^{f(t)}\: e^{i \alpha_-} du \right) ,
\eeq
and similarly for the approximate solution~$\tilde{U}_a(z)$.
We rewrite~\eqref{f} and~\eqref{ftt} as
\beql{fT}
f(t) - f(t_0) = T_1 + T_2 \:,\qquad
\tilde{f}(t) - f(t_0) = T_1 + T_3\:,
\eeq
where
\begin{align}
T_1 &:= \frac{1}{8} \left(8d+b\right) \tau^2 \\
T_2 &:= \frac{b}{4} \,\Big( \tau - \log(1+\tau) -\frac{\tau^2}{2} \Big) \label{T2def} \\
T_3 &:= - \frac{b}{12}\: \tau^3\:.
\end{align}

\begin{figure}
\begin{center}
\includegraphics[width=7cm]{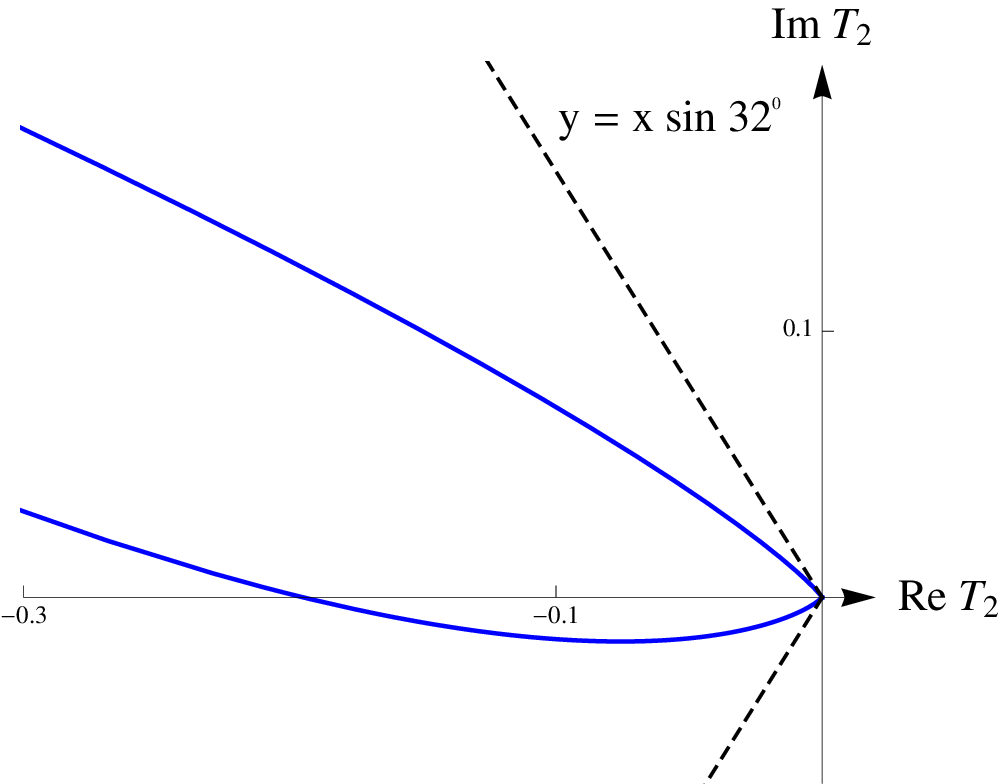} %$\quad$
\includegraphics[width=7cm]{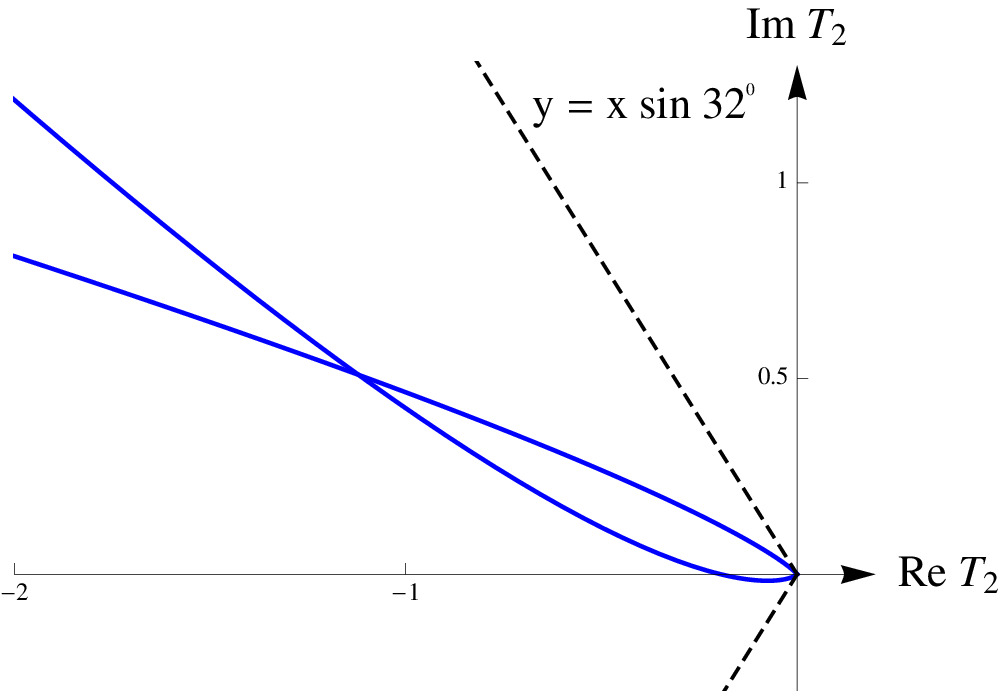} %$\quad$
\end{center}
\caption{The functions~$T_2^\pm$ in the complex plane.}
\label{figT21}
\end{figure}

\begin{figure}
\begin{center}
\includegraphics[width=7cm]{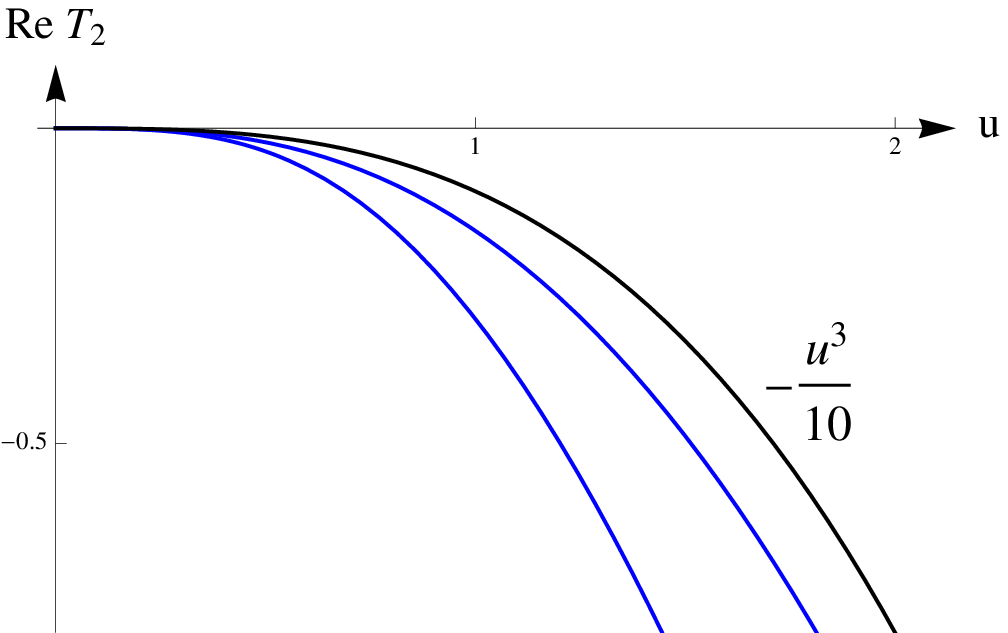}
\includegraphics[width=7cm]{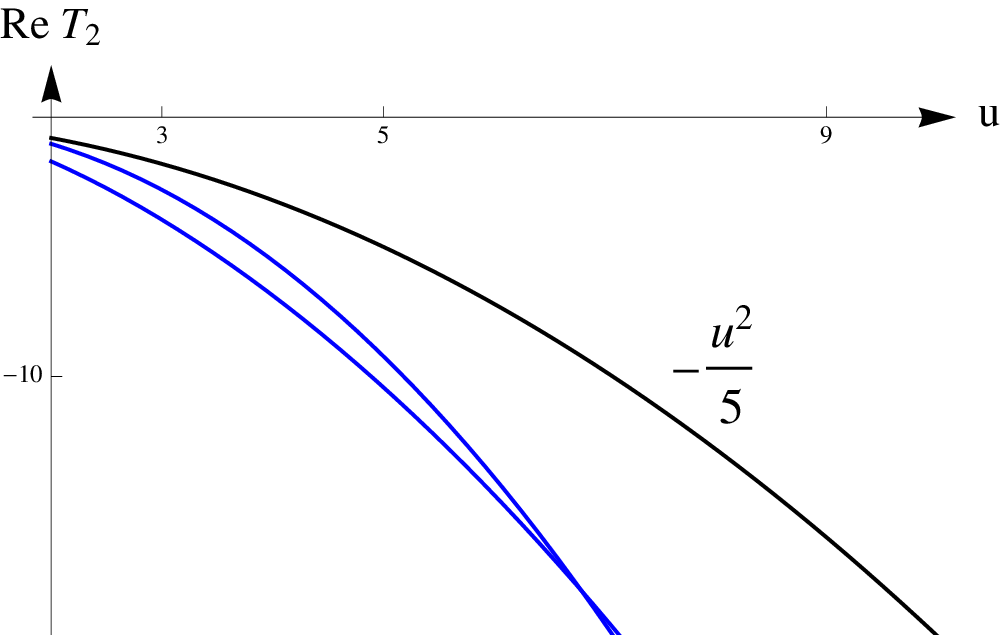}
\end{center}
\caption{Lower bounds for~$\re T_2^\pm$.}
\label{figT22}
\end{figure}

We now estimate~$T_1$, $T_2$ and~$T_3$. In order to estimate~$T_2$, we
use~\eqref{barg} to write~$b$ as
\[ b = i |b| \, e^{i \Delta} \qquad \text{with} \qquad |\Delta| \leq 2^\circ\:. \]
As in~\eqref{Umint}, we set~$\tau = e^{i \alpha_\pm} u$ with~$u\geq0$. It follows that
\beq \label{T1pmdef}
T_1 = e^{i \Delta} \: T_1^\pm \qquad \text{with} \qquad
T_1^\pm := \frac{1}{8}\:|b| \: e^{i (-90^\circ + 2 \alpha_\pm)} \left( \frac{4 t_0^2}{b} - 1 \right) u^2 \:.
\eeq
Using~\eqref{Rreg} and~\eqref{alphaval}, one sees that the points~$T_1^\pm$
lie in the regions~$\mathfrak{T}^\pm$ shown in Figure~\ref{figR}.
As a consequence, the points~$T_1^\pm$ always lie in the left half plane, and their angle
with the imaginary axis is at least~$5^\circ$. Since the factor~$e^{i \Delta}$ in~\eqref{T1pmdef}
describes a rotation
of at most~$2^\circ$, elementary trigonometry shows that
\beq \label{T1es}
\re T_1 \leq \frac{1}{2}\, \re T_1^\pm \leq -\frac{1}{2}\, \sin 5^\circ \: |T_1^\pm| 
\leq - \frac{1}{200} \left| 8d+b \right| u^2\:.
\eeq
In order to estimate~$T_2$, we first rewrite~\eqref{T2def} as
\[ T_2 = \frac{|b|}{4} \: e^{i \Delta}\: T_2^\pm \qquad \text{with} \qquad
T_2^\pm = i \Big( e^{i \alpha_\pm} \:u - \log(1 + e^{i \alpha_\pm} u) - e^{2 i \alpha_\pm}\: u^2 \Big) \,. \]
Using the explicit values of~$\alpha_\pm$ in~\eqref{alphaval}, one verifies that
the points~$T_2$ lie in the left half plane, and their angle with the imaginary axis is at least~$32^\circ$
(see Figure~\ref{figT21}).
Hence the angle of~$T_2$ to the imaginary axis is at least~$30^\circ$, and thus
\[ \re T_2 \leq \frac{|b|}{8}\: \re T_2^\pm\:. \]
Moreover, an explicit analysis (see the right of Figure~\ref{figT22}) shows that
\[ \re T_2^\pm \leq \left\{ \begin{array}{cl}
-u^3/10 & \text{if~$0 \leq u \leq 2$} \\[0.3em]
-u^2/5 & \text{if~$u > 2$\:.}  \end{array} \right. \]
We conclude that
\beql{T2es}
\re T_2 \leq -\frac{|b|}{80} \times \left\{ \begin{array}{cl}
u^3 & \text{if~$0 \leq u \leq 2$} \\[0.3em]
2 u^2 & \text{if~$u > 2$\:.}  \end{array} \right.
\eeq
To estimate~$T_3$, we note that
\[ T_3 = \frac{|b|}{12}\: e^{i (-90 + 3 \alpha_\pm + \Delta)}\: u^3 
=  \frac{|b|}{12} \times \left\{ \begin{array}{cl}
e^{i (234^\circ + \Delta)} & \text{in case~$+$} \\
e^{i (120^\circ + \Delta)} & \text{in case~$-$\:.}
\end{array} \right. \]
Thus~$T_3$ lies in the left half plane, and its angle with the imaginary axis is at least~$28^\circ$.
Hence
\beql{T3es}
\re T_3 \leq - \sin 28^\circ \: |T_3| \leq - \frac{|b|}{50}\: u^3\:.
\eeq
Using~\eqref{T1es}, \eqref{T2es} and~\eqref{T3es} in~\eqref{fT}, we conclude that
\begin{align}
\re \big( f(t) - f(t_0) \big) &\leq -\frac{1}{200} \left| 8d+b \right| u^2 
- \frac{|b|}{80} \times \left\{ \begin{array}{cl}
u^3 & \text{if~$0 \leq u \leq 2$} \label{fes} \\[0.3em]
2 u^2 & \text{if~$u > 2$\:.}  \end{array} \right. \\
\re \big( \tilde{f}(t) - f(t_0) \big) &\leq -\frac{1}{200} \left| 8d+b \right| u^2
-\frac{|b|}{50}\: u^3 \label{ftes}
\end{align}

As the estimate~\eqref{eb1} also holds in the case when~$t_0^2 \approx b/4$,
we do not want to make use of the first summands in~\eqref{fes} and~\eqref{ftes}.
For simplicity, we drop them,
\begin{align*}
\re \big( f(t) - f(t_0) \big) &\leq - \frac{|b|}{80} \times \left\{ \begin{array}{cl}
u^3 & \text{if~$0 \leq u \leq 2$} \\[0.3em]
2 u^2 & \text{if~$u > 2$\:.}  \end{array} \right. \\
\re \big( \tilde{f}(t) - f(t_0) \big) &\leq -\frac{|b|}{50}\: u^3 .
\end{align*}
Then
\begin{align*}
\int_L^\infty e^{\re \big( \tilde{f}(t) - f(t_0) \big)}\: du &\leq
\int_L^\infty e^{- \frac{|b|}{100}\: u^3}\: du \leq e^{- \frac{|b|\:L^3}{200}}
\int_0^\infty e^{-\frac{|b|}{200}\: u^3}\: du \\
&= e^{- \frac{|b|\, L^3}{200}} \:\Gamma
\Big( \frac{4}{3} \Big)\: \left( \frac{200}{|b|} \right)^\frac{1}{3}
\leq 6\: |b|^{-\frac{1}{3}}\: e^{- \frac{|b|\, L^3}{200}} \\
\int_L^\infty e^{\re \big( \tilde{f}(t) - f(t_0) \big)}\:du &\leq
\int_L^\infty \left( e^{- \frac{|b|}{80}\: u^3} + e^{- \frac{|b|}{40}\: u^2} \right)\: du \\
&\leq 6\: |b|^{-\frac{1}{3}}\: e^{- \frac{|b|\, L^3}{200}} +
e^{- \frac{|b|\, L^2}{80}} \int_0^\infty e^{- \frac{|b|\, L^2}{80}}\: du \\
&\leq 6\: |b|^{-\frac{1}{3}}\: e^{- \frac{|b|\, L^3}{200}} +
8\: |b|^{-\frac{1}{2}}\: e^{- \frac{|b|\, L^2}{80}} \:.
\end{align*}
On the interval~$[0,L]$, we can again follow the argument in~\eqref{smalles} to obtain
\[ \Big| e^{f(t)-f(t_0)} - e^{\tilde{f}(t)-f(t_0)} \Big| \leq |f(t) - \tilde{f}(t)|\:. \]
Using that~$\tilde{f}$ is the cubic Taylor approximation to~$f$, we obtain similar to~\eqref{quadratic} that
\[ |f(t) - \tilde{f}(t)| \leq \frac{|\tau|^4}{4!}\: \frac{3\, |b|}{2} \sup_{\tau} \frac{1}{|1+\tau|^4} 
\leq |b|\, u^4 \:. \]
Here in the last step we used that, according to~\eqref{alphaval}, on the contour~$\Gamma$
the inequality~$|1+\tau| \geq 1/2$ holds. As a consequence,
\beql{Lsmall}
\int_0^L \Big| e^{f(t)-f(t_0)} - e^{\tilde{f}(t)-f(t_0)} \Big| du \leq
\frac{|b|}{5}\, L^5 \:.
\eeq
Adding up all error terms, in view of~\eqref{Umint} we obtain for any~$L>0$ the inequality
\begin{align*}
\left| U_a^-(z) - \tilde{U}_a(z) \right| &\leq
2 \,e^{-\frac{\re z^2}{4} + f(t_0)} \:|t_0| 
\left( \frac{|b|}{5}\, L^5 + 12\: |b|^{-\frac{1}{3}}\: e^{- \frac{|b|\, L^3}{200}} +
8\: |b|^{-\frac{1}{2}}\: e^{- \frac{|b|\, L^2}{80}} \right)
\end{align*}
Choosing
\[ L = \frac{200}{|b|}\: \log \big(1+|b|^\frac{1}{3} \big) , \]
we obtain the estimate
\begin{align*}
\frac{e^{\frac{\re z^2}{4} - f(t_0)}}{|t_0|} \left| U_a^-(z) - \tilde{U}_a(z) \right| &\leq
1600 \cdot 5^{\frac{1}{3}}\: |b|^{-\frac{2}{3}}\: \log^\frac{5}{3} \big( 1+|b|^\frac{1}{3} \big)
+ 24 \, |b|^{-\frac{2}{3}} \\
&\qquad + 16\, |b|^{-\frac{1}{2}}\: \exp \Big(
-\frac{1}{4}\: 5^{\frac{1}{3}}\: |b|^\frac{1}{3}\: \log^\frac{2}{3} \big( 1+|b|^\frac{1}{3} \big) \Big) \\
&\leq 40 \:|b|^{-\frac{2}{3}} \left( 1 + 100\, \log^\frac{5}{3} \big( 1+|b|^\frac{1}{3} \big) \right) .
\end{align*}
This proves~\eqref{eb1}.

The estimate~\eqref{eb2} also applies in the case when~$|b|$ is small. This motivates us to drop the last
summands in~\eqref{fes} and~\eqref{ftes},
\[ \re \big( f(t) - f(t_0) \big),  \re \big( \tilde{f}(t) - f(t_0) \big)
\leq -\frac{1}{200} \left| 8d+b \right| u^2  \:. \]
Then
\begin{align*}
\int_L^\infty &e^{\re \big( f(t) - f(t_0) \big)}\, du \leq
\int_L^\infty e^{-\frac{|8d+b|}{200}\:  u^2} du \\
&\leq  e^{-\frac{|8d+b|}{400} \, L^2}
\int_0^\infty e^{-\frac{|8d+b|}{400}\,u^2} du
= \sqrt{2 \pi}\:5\, |8d+b|^{-\frac{1}{2}} \:e^{-\frac{|8d+b|}{400} \, L^2} \:,
\end{align*}
and similarly for~$\tilde{f}$. The integral over~$[0,L]$ can be estimated
again by~\eqref{Lsmall}. We thus obtain
\[ \left| U_a^-(z) - \tilde{U}_a(z) \right| \leq
2 \,e^{-\frac{\re z^2}{4} + f(t_0)} \:|t_0|  \left( 
\sqrt{2 \pi}\:10\, |8d+b|^{-\frac{1}{2}} \:e^{-\frac{|8d+b|}{400} \, L^2} +
\frac{|b|}{5}\, L^5 \right) . \]
Choosing
\[ L = 20\, |8d+b|^{-\frac{1}{2}}\: \log^\frac{1}{2} \Big( 1 + \frac{|8d+b|^2}{|b|} \Big) \:, \]
a short calculation gives~\eqref{eb2}.

\Thanks {{\em{Acknowledgments:}}
We would like to thank Martin Muldoon for several helpful remarks and for
introducing us to the literature of special functions.
We are grateful to the Vielberth Foundation, Regensburg, for generous support.

%\bibliographystyle{amsplain}
%\bibliography{../felix}

\begin{thebibliography}{10}

\bibitem{chandra}
S.~Chandrasekhar, \emph{The {M}athematical {T}heory of {B}lack {H}oles}, Oxford
  Classic Texts in the Physical Sciences, The Clarendon Press Oxford University
  Press, New York, 1998.

\bibitem{angular}
F.~Finster and H.~Schmid, \emph{Spectral estimates and non-selfadjoint
  perturbations of spheroidal wave operators}, J. Reine Angew. Math.
  \textbf{601} (2006), 71--107.

\bibitem{invariant}
F.~Finster and J.~Smoller, \emph{Error estimates for approximate solutions of
  the {R}iccati equation with real or complex potentials}, arXiv:0807.4406
  [math-ph], Arch. Ration. Mech. Anal. \textbf{197} (2010), no.~3, 985--1009.

\bibitem{tinvariant}
\bysame, \emph{Error estimates for approximate solutions of the angular
  {T}eukolsky equation in the {K}err geometry}, in preparation (2013).

\bibitem{tspectral}
\bysame, \emph{A spectral representation for spin-weighted spheroidal wave
  operators with complex aspherical parameter}, in preparation (2013).

\bibitem{hormanderI}
L.~H{\"o}rmander, \emph{The {A}nalysis of {L}inear {P}artial {D}ifferential
  {O}perators. {I}}, second ed., Grundlehren der Mathematischen Wissenschaften
  [Fundamental Principles of Mathematical Sciences], vol. 256, Springer-Verlag,
  Berlin, 1990, Distribution theory and Fourier analysis.

\bibitem{olver}
F.W.J. Olver, \emph{Asymptotics and {S}pecial {F}unctions}, Academic Press [A
  subsidiary of Harcourt Brace Jovanovich, Publishers], New York-London, 1974.

\bibitem{DLMF}
F.W.J. Olver, D.W. Lozier, R.F. Boisvert, and C.W. Clark (eds.), \emph{Digital
  {L}ibrary of {M}athematical {F}unctions}, National Institute of Standards and
  Technology from http://dlmf.nist.gov/ (release date 2011-07-01), Washington,
  DC, 2010.

\bibitem{temme+vidunas}
N.M. Temme and R.~Vidunas, \emph{Parabolic cylinder functions: {E}xamples of
  error bounds for asymptotic expansions}, arXiv:math/0205045v1 [math.CA],
  Anal. Appl. (Singap.) \textbf{1} (2003), no.~3, 265--288.

\bibitem{teukolsky}
S.A. Teukolsky, \emph{Perturbations of a rotating black hole {I}. {F}undamental
  equations for gravitational, electromagnetic, and neutrino-field
  perturbations}, Astrophys. J. \textbf{185} (1973), 635--647.

\end{thebibliography}
\providecommand{\bysame}{\leavevmode\hbox to3em{\hrulefill}\thinspace}
\providecommand{\MR}{\relax\ifhmode\unskip\space\fi MR }
% \MRhref is called by the amsart/book/proc definition of \MR.
\providecommand{\MRhref}[2]{%
  \href{http://www.ams.org/mathscinet-getitem?mr=#1}{#2}
}
\providecommand{\href}[2]{#2}

\end{document}